\documentclass[12pt]{article}

\usepackage[margin=0.55in]{geometry}
\usepackage{amsmath}
\usepackage{amssymb}
\usepackage{times}

\usepackage{bm}
\usepackage[round]{natbib}
\usepackage{color}
\usepackage{subcaption}
\usepackage{graphicx}
\usepackage{epstopdf}
\usepackage{chngcntr}
\usepackage{prodint}
\usepackage[onehalfspacing]{setspace}

\counterwithout{equation}{section}

\newtheorem{lemma}{Lemma}

\newtheorem{theorem}{Theorem}

\newcommand{\rnc}{\renewcommand}
\newcommand{\nc}{\newcommand}
\newcommand{\mrm}{\mathrm}

\nc{\black}{\color{black}}
\nc{\blue}{\color{black}}
\nc{\red}{\color{red}}
\nc{\mb}{\mathbb}
\nc{\mc}{\mathcal}
\nc{\E}{\mb{E}}
\nc{\N}{\mb{N}}
\nc{\R}{\mb{R}}
\nc{\Q}{\mb{Q}}
\rnc{\P}{\mrm P}
\rnc{\d}{\mrm d}
\nc{\C}{\mac{C}}
\nc{\D}{\mac{D}}
\nc{\B}{\mac{B}}
\nc{\oPo}{\stackrel{p}{\longrightarrow}}
\nc{\oWo}{\stackrel{w}{\longrightarrow}}
\nc{\oDo}{\stackrel{d}{\longrightarrow}}

\nc{\I}{\mc I}
\nc{\J}{\mc J}
\nc{\JI}{{\J\I}}
\rnc{\IJ}{{\I\J}}
\nc{\ji}{{\J|\I}}
\nc{\gDg}{\stackrel{d}{=}}
\nc{\nae}{Nelson-Aalen estimator}
\nc{\aje}{Aalen-Johansen estimator}
\nc{\naeL}{Nelson-Aalen estimator\ }
\nc{\ajeL}{Aalen-Johansen estimator\ }
\nc{\CIF}{cumulative incidence function}
\nc{\CIFL}{cumulative incidence function\ }
\nc{\wh}{\widehat}
\nc{\pr}{\textnormal{pr}}
\nc{\cred}{\color{red}}
\nc{\var}{\textnormal{var}}

\newcommand\blfootnote[1]{%
  \begingroup
  \renewcommand\thefootnote{}\footnote{#1}%
  \addtocounter{footnote}{-1}%
  \endgroup
} 

\begin{document}


\title{\vspace{-0cm}
{\blue Time-dynamic inference for non-Markov transition probabilities under independent right-censoring}}

\author{D. Dobler${}^*$
and A. C. Titman${}^\dagger$}


\maketitle

\allowdisplaybreaks

\begin{abstract} \blfootnote{${}^*${\blue Department of Mathematics, Faculty of Science, Vrije Universiteit Amsterdam, 1081 HV Amsterdam, Netherlands, e-mail:  d.dobler@vu.nl}}
\blfootnote{${}^\dagger$Department of Mathematics \& Statistics, Lancaster University, Bailrigg, Lancaster, LA1 4YW, United Kingdom, e-mail: a.titman@lancaster.ac.uk}
 \noindent In this article, weak convergence of the general non-Markov state transition probability estimator by \cite{titman15} is established
 which, up to now, {\blue has not been verified yet for other general non-Markov estimators.}
 A similar theorem is shown for the bootstrap, yielding resampling-based inference methods for statistical functionals.
 Formulas of the involved covariance functions are presented in detail.
 {\blue Particular applications include the conditional expected length of stay in a specific state, given occupation of another state in the past, as well as the construction of time-simultaneous confidence bands for the transition probabilities.}
 The expected lengths of stay in the two-sample liver cirrhosis data-set by \cite{abgk93} 
 are compared and confidence intervals for their difference are constructed.
 With borderline significance and in comparison to the placebo group, 
 the treatment group has an elevated expected length of stay in the healthy state given an earlier disease state occupation.
 In contrast, the Aalen--Johansen estimator-based confidence interval, which relies on a Markov assumption,
 leads to a drastically different conclusion.
 {\blue Also, graphical illustrations of confidence bands for the transition probabilities  demonstrate the biasedness of the Aalen-Johansen estimator in this data example.}
 The reliability of these results is assessed in a simulation study.
\end{abstract}


\section{Introduction}

Transition probabilities are essential quantities in survival analytic examinations of multi-state models.
Under independent right-censoring, for instance, 
the Aalen--Johansen estimator assesses these probabilities optimally in nonparametric Markovian multi-state models; cf. \cite{aalen78} and Section~IV.4 in \cite{abgk93} for more details.
The Markov assumption is crucial for the validity of the Aalen--Johansen estimator.
However, in real world applications, it may be unrealistic to postulate such a structure:
given the present patient's state, future developments are often not independent of a severe past illness history.
There have been several attempts to circumvent the Markov assumption:
State occupation probabilities in general models have been estimated by \cite{datta01}, \cite{datta02}, and \cite{glidden02} using the Aalen--Johansen estimator.
\cite{pepe1991b} and \cite{pepe1991} used a combination of two Kaplan--Meier estimators to assess the relevant transition probability in an illness-death model with recovery.
In illness-death models without recovery, \cite{meira-machado06} used Kaplan--Meier-based techniques and consecutive latent transition times to estimate transition probabilities.
More efficient variants of these estimators have been developed in \cite{deUna-alvarez15}.
A Kendall's $\tau$-based test of the Markov assumption in this progressive illness-death model was derived in \cite{rodriguez12}.
A competing risks-based estimator was proposed by \cite{allignol14} which relies on strict censoring assumptions.
Eventually, \cite{titman15} and \cite{putter16n} developed transition probability estimators in general non-Markov models by using different transition probability decompositions and utilizations of the Aalen--Johansen estimator.
However, weak convergence properties of these estimators as elements of c\`adl\`ag function spaces have not been analyzed yet.

In the present paper, we focus on the non-Markov state transition probability estimator by \cite{titman15}
and its weak convergence for the following reasons:
the estimator has a simple, but intuitive structure
and it allows estimation of general transition probabilities between sets of states rather than single states.
Finally, the bootstrap is shown to correctly reproduce the weak limit distribution.
Throughout the article, let $(\Omega, \mathcal A, P)$ be the underlying probability space.
To introduce the \cite{titman15} estimator,
we consider an independently right-censored multi-state model with $R+1 \in \N$ different states.
That is, the data consists of independent copies $X_1, \dots, X_n$ of a multi-state process which are indexed by time $t \geq s$ and which occupy the states $0,1,\dots, R$.
The right-censoring is modeled by independent random variables $C_1, \dots, C_n$ after which further observation of the corresponding multi-state processes is not possible any more.
  Let $\I$ and $\J$ be subsets of $\{0,1, \dots, R\}$ for which the transition probabilities 
  $$ P_{\I\J}(s,t) = P(X(t) \in \J  \mid  X(s) \in \I)$$
  shall be estimated.
  Here and below, irrelevant subscripts are removed for notational convenience.
  Introduce the set of states $\mc A_\J \subset \J$, that implies sure occupancy of $\J$ at all later points in time,
  and the set of states $\mc R_\J \subset \{0,1,\dots, R\} \setminus \J$ which prevents occupancy of $\J$ at all later points of time.
  Individual competing risks processes, $Z_i, i=1,\dots,n$, are used in an intermediate estimation step.
  $Z_i(t) = 0$ 
  indicates whether the multi-state process $i$ is at time $t$ in none of the absorbing subsets $\mc A_\J$ or $\mc R_\J$.
  On the other hand, $Z_i(t) = 1$ holds if the transition to $\mc A_\J$ has been observed until time $t$
  and similarly $Z_i(t) = 2$ for the transition to $\mc R_J$.
  Let $Y_i(t), i=1,\dots,n$ denote the at-risk indicators of the competing risks processes
  and let $\delta_i, i=1,\dots,n$ be the corresponding censoring indicators, 
  i.e., $\delta_i = 1$ if the transition of $Z_i$ is eventually observed and $\delta_i = 0$ if a censoring comes first.
  The final data-set is $\chi = \{ X_i, Z_i, Y_i, \delta_i: i=1,\dots,n \}$
  whereof, in case of $\delta_i=0$, the processes $X_i$ and $Z_i$ are observable only until censoring occurs.
  Hence, incorporating random right-censoring, the transition probability is decomposed into 
  \begin{align*}
    P_{\I\J}(s,t)  = & P( Z(t) = 1 \mid Y(s) = 1, X(s) \in \I) \\
    & + P( Y(t)=1 \mid Y(s) = 1, X(s) \in \I ) \\
    & \quad \times P( X(t) \in \J \mid Y(t) = 1, X(s) \in \I ) \\
    = & F_1(t) + F_0(t) p_\ji(t).
  \end{align*}
  The first two components in this decomposition are estimated by the standard Kaplan--Meier estimator $\wh F_0$ for the survival function and the Aalen--Johansen estimator $\wh F_1$ for the cumulative incidence function of the sub-group of those individuals which were observed to occupy state $\I$ at time $s$.
  The remaining conditional probability is estimated using the empirical proportion
  $$ \wh p_{\mc J \mid \mc I}(t) = \frac{\wh p_{\mc J \mc I}(t)}{\wh p_{\mc I}(t)} = \frac{n^{-1} \sum_{\ell=1}^n 1 \{ X_\ell(s) \in \I,  X_\ell(t) \in \J, Y_\ell(t) = 1 \}}
   {n^{-1} \sum_{\ell=1}^n 1 \{ X_\ell(s) \in \I,  Y_\ell(t) = 1 \}} ,$$
   where $1 \{\cdot\}$ denotes the indicator function.
   Therefore, the overall transition probability estimator by \cite{titman15} is given as
   \begin{align*}
    \wh P_{\I\J}(s,t) = \wh F_1(t) + \wh F_0(t) \wh p_\ji(t).
   \end{align*}
  This article shall establish the weak convergence of $\sqrt{n} ( \wh P_{\I\J}(s, \cdot) - P_{\I\J}(s, \cdot) )$ and of a bootstrap variant thereof as $n \rightarrow \infty$.
  Finally, implications to confidence interval construction for the conditional expected length of stay
  {\blue and to time-simultaneous confidence bands for the transition probabilities}
  are demonstrated.
  {\blue The proofs of all Theorems and Lemmata as well as detailed derivations of all asymptotic covariance functions are given in the appendix.}

\section{Main Results}
\label{sec:main}

We focus on the \cite{titman15} estimator during $[s,\tau]$
where $\tau \in (s,\infty)$ is a terminal time satisfying $F_0(\tau) > 0$ and $P(C > \tau) > 0$.
The large sample properties of $\wh P_{\I\J}(s, \cdot)$ are derived from central limit theorems of its individual components.
To obtain such for the fraction process $\wh p_\ji$, we throughout assume the existence of bounded one- and two-step hazard functions for {\blue transitions into} $\J$ {\blue and $\bar\J = \{0,1,\dots,R\} \setminus \J$, respectively,}
which are approached uniformly:
  \begin{align}
   \lim_{\delta \downarrow 0} \sup_{t-t_1 = \delta} \Big| & \frac{P (X(t) \in \J  \mid  X(t_1) \notin \J )}{t - t_1} 
    - {\blue \lambda_{\bar{\J}\J}(t)} \Big| = 0, \label{eq:main_cond_1} \\
   \lim_{\delta \downarrow 0} \sup_{t-t_1 = \delta} \Big| & \frac{P (X(t) \notin \J \mid X(t_1) \in \J )}{t - t_1} 
    - {\blue \lambda_{\J\bar{\J}}(t)} \Big| = 0, \\ 
   \lim_{\delta \downarrow 0} \sup_{\substack{t_2-t = \delta \\ t_1 < t}} \Big| & \frac{P ( X(t_2) \in \J \mid X(t) \notin \J, X(t_1) \in \J )}{t_2 - t} 
    - {\blue\lambda_{\bar{\J}\J\mid \J}(t_2 \mid t_1)} \Big| = 0, \\
   \lim_{\delta \downarrow 0} \sup_{\substack{t_2-t = \delta \\ t_1 < t}} \Big| & \frac{P(X(t_2) \notin \J \mid X(t) \in \J, X(t_1) \notin \J )}{t_2 - t} 
    - {\blue\lambda_{\J\bar{\J}\mid\bar{\J}}(t_2 \mid t_1)} \Big| = 0. \label{eq:main_cond_4}
  \end{align}
  Here, $t_1, t, t_2 \in (s,\tau]$.
  In contrast to the standard choice in the literature, we let the at-risk indicators $Y_i,i=1, \dots, n$,
  that are involved in the estimator $\wh p_\ji$,
  be right-continuous.
  The typical reason for using the left-continuous version is the applicability of martingale methods
  which do not matter in our analysis of the estimators.
  On top of that, this choice does not matter in the present article as all transition times are continuously distributed.
  The definitions of the Kaplan--Meier and the Aalen--Johansen estimators remain unchanged.
  We equip the c\`adl\`ag function space $D[s,\tau]$ with the supremum norm
  and products thereof with the maximum-supremum norm;
  see Chapter~3.9 in \cite{vaart96} for applications in survival analysis with this choice.
  In order to state our main theorem, it is convenient to provide an auxiliary lemma:
  \begin{lemma}
   \label{lemma:main}
   Under conditions \eqref{eq:main_cond_1} to \eqref{eq:main_cond_4} and as $n \rightarrow \infty$,
   the following joint convergence holds in distribution on $D^3[s,\tau]$:
   \begin{align*}
     \sqrt{n} (\wh F_0 - F_0, \wh F_1 - F_1, \wh p_\ji - p_\ji) \rightarrow (L_0, L_1, G_\ji)
   \end{align*}
   Here, $L_0$, $L_1$, and $G_\ji$ are zero-mean Gaussian processes. 
   Denote their variance-covariance functions, arranged in a convenient $(3 \times 3)$-matrix, by
   \begin{align*}
    \begin{pmatrix}
     \Sigma_{00} & \Sigma_{01} & \Sigma_{0G} \\
     \Sigma_{10} & \Sigma_{11} & \Sigma_{1G} \\
     \Sigma_{G0} & \Sigma_{G1} & \Sigma_{GG}
    \end{pmatrix} : [s, \tau]^2 \rightarrow \R^{3 \times 3}.
   \end{align*}
\end{lemma}
    {\blue Detailed formulas for these covariance functions are presented in the appendix.
    From the presentation therein it is apparent that the covariances can be estimated straightforwardly.}

Define $r \wedge t = \min(r,t)$ and $r \vee t = \max(r,t)$.
An application of the functional delta-method immediately yields our main theorem:
\begin{theorem}
 \label{thm:main}
 Under conditions \eqref{eq:main_cond_1} to \eqref{eq:main_cond_4} and as $n \rightarrow \infty$,
 the following convergence holds in distribution on $D[s,\tau]$:
 \begin{align}
 \label{eq:main}
  \sqrt{n} ( \wh P_{\I\J}(s, \cdot ) - P_{\I\J}(s,\cdot) ) \rightarrow U
 \end{align}
 where $U$ is a zero-mean Gaussian process with covariance function $\Gamma_{\I\J}$ given by
 \begin{align*}
  (r,t) & \mapsto p_\ji(r) p_\ji(t) \Sigma_{00}(r,t) + \Sigma_{11}(r,t) + F_0(r) F_0(t) \Sigma_{GG}(r,t) \\
    & + p_\ji(r)\Sigma_{01}(r,t) 
    + p_\ji(t) \Sigma_{10}(r,t) \\
    & + p_\ji(r \vee t) F_0(r \wedge t) \Sigma_{0G}(r \vee t, r \wedge t) 
    + F_0(r \wedge t) \Sigma_{1G}(r \vee t,r \wedge t) .
 \end{align*}
\end{theorem}
To be able to use this result for time-simultaneous inference on $P_{\I\J}(s,\cdot)$ 
or simply for improving the accurateness of $P_{\I\J}(s,\cdot)$-based inference methods with a pivotal limit distribution,
the classical bootstrap seems to be a useful choice of a resampling method.
To this end, a bootstrap sample $\chi^*$ is obtained by {\blue randomly drawing quadruples $(X_i, Z_i, Y_i,\delta_i)$ $n$ times} with replacement from the original sample $\chi$.
This bootstrap sample is used to recalculate the transition probability estimator, resulting in
\begin{align*}
 P_{\I\J}^*(s,t) = F_1^*(t) + F_0^*(t) p_\ji^*(t), \quad t \in [s,\tau].
\end{align*}
For this bootstrap counterpart, a central limit theorem similar to the main Theorem~\ref{thm:main} holds:
\begin{theorem}
\label{thm:bs}
 Under conditions \eqref{eq:main_cond_1} to \eqref{eq:main_cond_4}, as $n \rightarrow \infty$,
 and given $\chi$,
 the following conditional convergence in distribution on $D[s,\tau]$ holds in probability:
 \begin{align}
 \label{eq:bs}
  \sqrt{n} ( P^*_{\I\J}(s, \cdot ) - \wh P_{\I\J}(s,\cdot) ) \rightarrow U
 \end{align}
 where $U$ is the same Gaussian process as in Theorem~\ref{thm:main}.
\end{theorem}
That is, the bootstrap succeeds in correctly recovering the limit distribution of the transition probability estimate $\wh P_{\I\J}(s,\cdot)$, enabling an approximation of its finite sample distribution.

{\blue 
\section{Inference Procedures}
\label{sec:elos}
\subsection{Conditional Expected Length of Stay}
}

  A practical application of both main theorems is a two-sample comparison of expected lengths of stay
  in the state set of interest, $\J$, given occupation of a fixed state set $\I$ at time $s \in [0,\tau]$.
  Considering temporarily the above one-sample set-up, the conditional expected length of stay is
  \begin{align*}
   e_{\IJ}(s, \tau) & = E \Big[ \int_s^\tau 1\{ X(u) \in \J \} {\color{blue}\d u} \ | \ X(s) \in \I  \Big] \\
    & = \int_s^\tau P ( X(u) \in \J \ | \ X(s) \in \I ) \d u;
  \end{align*}
  see e.g. \cite{grand16} for a pseudo-observation regression treatment.
  {\color{black}
  An R package for deriving the change in lengths of stay, that relies on the Aalen--Johansen estimator, is described in \cite{wangler06}.
  For further literature on lengths of stay within multi-state models, see the articles cited therein.
  In the special case of a simple survival model, the only expected length of stay of interest is the mean residual lifetime; see e.g. \cite{meilijson72}.
  Here, however, the Markov assumption shall not be imposed such that estimation of the length of stay is based on the \cite{titman15} estimator.}
  
  Sample-specific quantities are furnished with a superscript ${}^{(\ell)},\ell = 1,2$.
  For instance, $\wh P^{(1)}_{\IJ}(s,t)$ is the first sample's transition probability estimate 
  to switch into some state of $\J$ at time $t$ given that a state of $\I$ is occupied at time $s$.
  The sample sizes $n_1$ and $n_2$ may be different.
  We would like to test the null hypothesis
  $$ H_{\J \mid \I} : {e}^{(1)}_{\I\J}(s,\tau) \geq {e}^{(2)}_{\I\J}(s,\tau) \quad  \text{versus} \quad 
    K_{\J \mid \I} : {e}^{(1)}_{\I\J}(s,\tau) < {e}^{(2)}_{\I\J}(s,\tau). $$
    {\blue If $\J$ is a favourable set of states, then the rejection of $H_{\J \mid \I}$ is an indication of the} superiority of the second treatment over the first.
    A reasonable estimator {\blue of} ${e}^{(\ell)}_{\I\J}(s,\tau)$ is
    $
     \wh e_{\I\J}^{(\ell)}(s,\tau) = \int_s^\tau \wh P_{\I\J}^{(\ell)}(s,u) \d u,
    $
    \ $\ell =1,2,$.
  Note that, even under the boundary $\partial H_{\J|\I} :  {e}^{(1)}_{\IJ}(s,\tau) = {e}^{(2)}_{\IJ}(s,\tau)$, 
  the underlying transition probabilities $P^{(1)}_{\I\J}$ and  $P^{(2)}_{\I\J}$ are allowed to differ.
  {\blue Hence, we do not rely on the restrictive null hypothesis $\tilde H: P^{(1)}_{\I\J} \equiv P^{(2)}_{\I\J}$
  but rather on the hypothesis $H_{\J \mid \I}$ of actual interest.}
  \begin{theorem}
   \label{thm:elos}
   Under conditions \eqref{eq:main_cond_1} to \eqref{eq:main_cond_4} and if \ $0 < \liminf n_1 / n_2 \leq \limsup n_1/n_2$ $< \infty$  as $n \rightarrow \infty$, the following convergence in distribution holds:
   \begin{align*}
    \sqrt{\frac{n_1 n_2}{n_1 + n_2}} \ \frac{\wh {e}_{\IJ}^{(1)}(s,\tau) - \wh {e}_{\IJ}^{(2)}(s,\tau)}{\wh \sigma_{\I\J}} \rightarrow N(0,1). 
   \end{align*}
   Here, $\wh \sigma^2_{\I\J}$ is any consistent estimator, e.g., found via plug-in, of the asymptotic variance
   $$ \sigma^2_{\I\J} = (1 - \lambda) \int_s^\tau \int_s^\tau \Gamma_\IJ^{(1)}(u,v) \d u \d v 
      + \lambda \int_s^\tau \int_s^\tau \Gamma_\IJ^{(2)}(u,v) \d u \d v $$
   along subsequences of $n_1 / (n_1 + n_2 ) \rightarrow \lambda \in (0,1)$.
  \end{theorem}
  Replacing all estimators in the previous theorem by their bootstrap counterparts,
  where both underlying bootstrap samples $\chi^{(1)*}$ and $\chi^{(2)*}$ are obtained
  by $n_1$ and $n_2$ times independently and separately drawing {\blue with replacement} from $\chi^{(1)}$ and $\chi^{(2)}$, respectively,
  a similar convergence in conditional distribution given $\chi^{(1)} \cup \chi^{(2)}$ holds in probability.   
  This result may be used for bootstrap-based and asymptotically exact confidence intervals for the difference in the conditional expected lengths of stay, $ {e}_{\IJ}^{(1)}(s,\tau) - {e}_{\IJ}^{(2)}(s,\tau) $;
  see Sections~\ref{sec:simus} and~\ref{sec:example} below.
  
  {\blue
  \subsection{Simultaneous Confidence Bands}
  \label{sec:scb}
  Even though the bootstrap procedure is useful but not strictly necessary for inference on the expected length of stay, 
  it plays an essential role in time-simultaneous inference on the transition probabilities $t \mapsto P_{\I\J}(s,t)$.
  This is due to the unknown stochastic behaviour of the limit process $U$ which, again in the context of creating confidence bands, is also the reason for the inevitableness of resampling procedures for Aalen-Johansen estimators even if the Markov assumption is fulfilled; see \cite{bluhmki18} for theoretical justifications and their practical performance.
  
  For the derivation of reasonable time-simultaneous $1-\alpha \in (0,1)$ confidence bands for  $P_{\I\J}(s,\cdot)$ on a time interval $[t_1,t_2] \subset [s,\tau]$ based on  $\wh P_{\I\J}(s,\cdot)$ in combination with the bootstrap,
  we focus on the process 
  $$ B_n : t \longmapsto \sqrt{n} \cdot  w(t) \cdot (\phi(\wh P_{\I\J}(s,t)) - \phi(P_{\I\J}(s,t))), \quad t \in [t_1,t_2], $$
  where a transformation $\phi$ will ensure bands within the probability bounds $0\%$ and $100\%$,
  and a suitable multiplicative weight function $w$ will stabilize the bands for small samples.
  In particular, we use one of the log-log transformations $\phi_1(p) = \log(-\log(p))$ or $\phi_2(p) = \log(-\log(1-p))$ \citep{lin97}, depending on whether we expect $P_{\I\J}(s,t_1)$ to be closer to 1 or 0.
  As weight function, for transformation $k\in\{1,2\}$, we choose $w_{1k}(t) = [\phi_k'(\wh P_{\I\J}(s,t)) \cdot \wh\sigma_\IJ(t)]^{-1}$,
  \begin{align*}
   w_{21}(t) & = -\Big(\phi_1'(\wh P_{\I\J}(s,t)) \cdot \wh P_{\I\J}(s,t) \cdot \Big[1 + \frac{\wh\sigma^2_\IJ(t)}{\wh P_{\I\J}^2(s,t)}\Big] \Big)^{-1} \\
  & = -\log(\wh P_{\I\J}(s,t)) / \Big[1 + \frac{\wh\sigma^2_\IJ(t)}{\wh P_{\I\J}^2(s,t)}\Big]^{-1}, \\
\text{or} \quad    w_{22}(t) & = \Big(\phi_2'(\wh P_{\I\J}(s,t)) \cdot  (1-\wh P_{\I\J}(s,t)) \cdot \Big[1 + \frac{\wh\sigma^2_\IJ(t)}{(1-\wh P_{\I\J}(s,t))^2}\Big] \Big)^{-1} \\
    & =  -\log(1 - \wh P_{\I\J}(s,t)) /   \Big[1 + \frac{\wh\sigma^2_\IJ(t)}{(1-\wh P_{\I\J}(s,t))^2}\Big]^{-1}.
  \end{align*}
  With these choices, the resulting confidence bands will correspond to the classical equal precision (EP) bands for $w=w_{1k}$ or the Hall-Wellner (H-W) bands for $w=w_{2k}$ if $\J$ is an absorbing subset of states;
  cf. p. 265ff. in \cite{abgk93} for the survival case and $k=1$ or \cite{lin97} in the presence of competing risks and $k=2$.
  
  The confidence bands are found by solving for $P_{\I\J}(s,t)$ in the probability
  $P ( \sup_{t \in [t_1,t_2]} | B_n(t)| \leq q_{1-\alpha} ) = 1-\alpha $.
  Here, the value of $q_{1-\alpha}$ is approximated by the $(1-\alpha)$-quantile $q_{1-\alpha}^*$ of the conditional distribution of the supremum of the bootstrap version $B_n^*(t) = \sqrt{n} \cdot  w^*(t) \cdot (\phi(P_{\I\J}^*(s,t)) - \phi( \wh P_{\I\J}(s,t)))$ of $B_n$ given the data, 
  where in the definition of $w^*$ all estimators $\wh P_{\I\J}(s,t)$ and $\wh\sigma_\IJ(t)$ are replaced with their bootstrap counterparts $P_{\I\J}^*(s,t)$ and $\sigma^{*}_\IJ(t)$.
  
  Finally, the resulting $1-\alpha$ confidence bands for $t \mapsto P_{\I\J}(s,t)$ are as follows:
  \begin{align*}
   CB_{1k} & = [ \phi_1^{-1} ( \phi_1(\wh P_{\I\J}(s,t)) \mp q_{1-\alpha}^* / (\sqrt{n}  w_{1k}(t)) ) ]_{t_1 \leq t \leq t_2} \\
   & = [ \wh P_{\I\J}(s,t)^{\exp(\mp q_{1-\alpha}^* / (\sqrt{n}  w_{1k}(t)) )} ]_{t_1 \leq t \leq t_2} \\
   \text{and}\quad 
   CB_{2k} & = [ \phi_2^{-1} ( \phi_2(\wh P_{\I\J}(s,t)) \mp q_{1-\alpha}^* / (\sqrt{n}  w_{2k}(t)) ) ]_{t_1 \leq t \leq t_2} \\
   & = [ 1 - (1 -  \wh P_{\I\J}(s,t))^{\exp(\mp q_{1-\alpha}^* / (\sqrt{n}  w_{2k}(t)) )} ]_{t_1 \leq t \leq t_2}
  \end{align*}
  The performance of these confidence bands as well as the confidence band obtained from not transforming or weighting the transition probability estimator
  are assessed in a simulation study in the subsequent sections.
  }


\section{Simulation Study}
\label{sec:simus}
{\blue
\subsection{Conditional Expected Length of Stay}
}
To assess the performance of the conditional length of stay estimators for finite sample sizes a small simulation study is conducted. Data are simulated from a pathological non-Markov model in which the subsequent dynamics of the process depend on the state occupied at $t=4$. A three-state illness-death model with recovery is used
{\blue where the state $0$ means ``healthy", $1$ ``ill", and $2$ ``dead".
The} transition intensities are $\lambda_{02}(t) = 0.02, \lambda_{10}(t) = 0.3, \lambda_{12}(t)=0.1$, while
\begin{equation}\lambda_{01}(t) =\begin{cases} 0.3 &\mbox{if}~X(4)=1, t>4\\ 0.6 & \mbox{otherwise.} \end{cases}
\end{equation}
{\blue In the last part of the appendix we prove that the conditions~\eqref{eq:main_cond_1}--\eqref{eq:main_cond_4} are satisfied in the present set-up.}
Subjects are independently right-censored via an exponential distribution with rate 0.04. 
{\blue We focus on the expected length of stay in the health state conditional on an earlier illness at time $s=5$, i.e. $\I = \{1\}$ and $\J = \{0\}$.}
For each simulated dataset, $e_{10}(5,30)$ is estimated and 95\% confidence intervals are constructed via three methods; a Wald interval using the plug-in variance estimator $\widehat{\sigma}_{10}$ {\blue of ${\sigma}_{10}$, in which the canonical estimators of $p_{0 \mid 1}$, $F_0$, $\Sigma_{00}$, $\Sigma_{0G}$, and $\Sigma_{GG}$ are chosen,} a classical bootstrap and a bootstrap-t procedure using {\blue the bootstrap-version of} $\wh{\sigma}_{10}$ to studentize the bootstrap samples. For each bootstrap, $B=1000$ samples are generated. 
%
%
%

%
%
%
In the present {\blue three-state model}, the \cite{titman15} transition probability estimator reduces to ${\blue \wh P_{10}(5,t)=} \wh P_{\I\J}(5,t) = \wh F_0(t) \wh p_\ji(t) {\blue = \wh F_0(t) \wh p_{0\mid1}(t)}$ {\blue because the healthy state is non-absorbing}.
Thus, the asymptotic variance of each sample-specific conditional expected length of stay estimator is
\begin{align*}
 \int_s^\tau \int_s^\tau \Big( & p_{0 \mid 1}(u) p_{0 \mid 1}(v) \Sigma_{00}(u,v) + F_0(u) F_0(v) \Sigma_{GG}(u,v) \\
  & + p_{0 \mid 1}(u \vee v) F_0(u \wedge v) \Sigma_{0G}(u \wedge v, u \vee v) \Big) \d u \d v.
\end{align*}
As a competing method, Aalen--Johansen estimator-based Wald-type confidence intervals are constructed for the same expected length of stay.
Standard arguments yield that the asymptotic covariance function of the normalized Aalen--Johansen estimator 
for the $1 \rightarrow 0$ transition is $(r,t) \mapsto $
\begin{align*}
 \sum_{j=0}^2 \sum_{k \neq j} \int_s^{r \wedge t} P_{1j}^2(s,u) ( P_{k0}(u,t) - P_{j0}(u,t) ) 
 ( P_{k0}(u,r) - P_{j0}(u,r) ) \frac{\alpha_{jk}(u)}{y_j(u)} \d u;
\end{align*}
cf. Section~IV.4.2 in \cite{abgk93} for the corresponding particular variance formula.
The expected length of stay asymptotic variance again results from a double integral of this covariance function with integration range $r,t \in (s,\tau]$.

We consider several different scenarios for the total sample size, $n=50, 100,$ $150,$ and $200$. Note that the non-Markov estimator only uses the subgroup of subjects satisfying $C_i > 5, X_i(t)=1$ which corresponds to only $0.367n$ subjects, on average. For each of the scenarios, 10000 datasets are generated. As a result the empirical coverage percentages have approximate Monte-Carlo standard error of around 0.2.   

\begin{table}
\def~{\hphantom{0}}
\caption{Empirical bias and coverage of nominal 95\% confidence intervals for simulated datasets}%
{\blue
{
\centering
\begin{tabular}{lcccccc}
 &\multicolumn{2}{c}{Bias}&\multicolumn{4}{c}{Coverage \%}\\
$n$&AJ & NM& Wald (AJ) & Wald (NM) & Naive bootstrap & Bootstrap-t\\[5pt]
50 & -0.231 & 0.181 & 92.09 & 91.01 & 94.23 & 90.85\\
100& -0.268 & 0.085 & 89.95 & 94.54 & 94.08 & 94.43\\
150& -0.279 & 0.063 & 88.14 & 95.24 & 94.21 & 95.00\\
200& -0.297 & 0.033 & 86.34 & 95.56 & 94.54 & 95.43\\
\end{tabular}
}
}
\label{table_simulate}
\end{table}

The simulation results are shown in Table \ref{table_simulate}. Both the Aalen-Johansen and non-Markov estimator have substantial bias when $n=50$ leading to under coverage of all the nominal 95\% confidence intervals. For larger $n$, both of the bootstrap confidence intervals and the Wald-based interval for the non-Markov estimator give close to nominal coverage. The coverage of the Aalen-Johansen based Wald confidence interval deteriorates with increasing $n$ due to the estimator being inherently biased.

{\blue
\subsection{Simultaneous Confidence Bands}
\label{sec:bands}

To investigate the coverage probabilities of simultaneous confidence bands introduced in Section \ref{sec:scb}, datasets are simulated under the same non-Markov model as above, again considering total sample sizes of $n=50, 100, 150,$ and $200$ patients. 95\% simultaneous confidence bands are constructed for $p_{21}(5,t)$ for $t \in (6,7].$ For the non-Markov estimator, in addition to the Hall-Wellner and equal precision bands using the transformation $\phi_2(p)$, a naive confidence band based on a constant weight function, $w(t) \equiv 1$, and an identity transformation, $\phi(p) = p$, is also constructed. In addition, for the Aalen-Johansen estimates, EP bands based on $\phi_2(p)$, are constructed via the wild bootstrap using the R code which has been made available in the supplement to \cite{bluhmki18}. In all cases, for each bootstrap $B=1000$ samples are generated.

Table \ref{table_simulate2} gives the empirical coverage probabilities of the simultaneous confidence bands based on 5000 simulated datasets. For the bands based on the Non-Markov estimator, the H-W and EP bands tend to over-cover quite markedly for small sample sizes, while the naive bands under-cover. Adequate coverage is achieved for the H-W and EP bands by $n=200$. Since the Aalen-Johansen estimator is biased for this scenario we would expect the coverage of the EP Aalen-Johansen bands to deteriorate with increasing sample size. However, it appears this is counteracted by the estimated quantiles, $q^{*}_{1-\alpha}$ increasing with $n$. As a consequence, the coverage initially grows with increasing $n$ before beginning to decrease.

}

\begin{table}
\def~{\hphantom{0}}
\caption{Empirical coverage of nominal 95\% simultaneous confidence bands for $p_{21}(5,t)$ for $t \in (6,7]$. H-W = Hall-Wellner, EP = Equal precision}%
{\centering
{\blue


\begin{tabular}{lcccc}
 &\multicolumn{3}{c}{Non-Markov}&Aalen-Johansen\\
$n$&H-W & EP& Naive & EP\\[5pt]
50 & 0.9996 & 0.9972 & 0.9226 & 0.9154 \\
100& 0.9896 & 0.9846 & 0.9238 & 0.9188\\
150& 0.9604 & 0.9542 & 0.9298 & 0.9230  \\
200& 0.9496 & 0.9496 & 0.9420 & 0.9124  \\
\end{tabular}
}

}
\label{table_simulate2}
\end{table}

\section{Application to the Liver Cirrhosis Data-Set}
\label{sec:example}

As an illustrative example we consider the liver cirrhosis data set introduced by \cite{abgk93} and also analyzed by \cite{titman15}. Patients were randomized to either a treatment of prednisone (251 patients) or a placebo
(237 patients). A three-state illness-death model with recovery is assumed where the healthy and illness states correspond to normal and elevated levels of prothrombin, respectively. A potential measure of the effectiveness of prednisone as a treatment is the expected length-of-stay in the normal prothrombin level state, from a defined starting point.
Specifically we consider
$\Delta e_{10} = e_{10}^{(1)}(s,\tau) - e_{10}^{(2)}(s,\tau)$, the difference in conditional expected length-of-stay in normal prothrombin levels, given the patient is alive and with abnormal prothrombin at time $s$, for the prednisone and placebo groups, corresponding to $\ell=1$ and $\ell=2$, respectively. 
\cite{titman15} observed an apparent treatment difference with respect to the transition probabilities with respect to ${\blue s}=1000$ days post-randomization. 
We consider the expected length-of-stay in state 0 up to time $\tau = 3000$ days. 

In this case, $n_1 = 26$ patients in the prednisone group and $n_2 = 35$ patients in the placebo group meet the condition that $X_i(1000)=1$ and $Y_i(1000)=1$.  Bootstrap based confidence intervals are constructed using $1000$ bootstrap samples within each group.

\begin{table}
\def~{\hphantom{0}}
\caption{Estimated difference in conditional length-of-stay for liver cirrhosis data}%
\begin{center}
\begin{tabular}{lcc}
Method& $\Delta e_{10}$&95\% CI\\[5pt]
Wald & 375.3 &(-2.1, 752.8)\\
Naive bootstrap&375.3 &(14.7, 740.0)\\
Bootstrap-t &375.3 & (-5.1, 747.3)\\[3pt]
AJ (Wald)& 11.9 & (-217.3, 241.1)\\
\end{tabular}
\end{center}
\label{tableliverelos}
\end{table}

Table \ref{tableliverelos} shows the three constructed 95\% confidence intervals, which are broad\-ly similar, although the naive bootstrap interval excludes 0, whereas the Wald and bootstrap-t intervals do not. The Aalen-Johansen based estimate and a Wald confidence interval is also given. It is seen that there is a dramatic difference between the non-Markov and Aalen-Johansen estimates, with the latter indicating no treatment difference. Potentially, the apparent non-Markov behaviour in the data may be due to patient heterogeneity within the treatment groups.


{\blue
To illustrate the construction of simultaneous confidence bands we construct 95\% confidence bands for $P_{21}(500,t)$, in the interval $t \in (750,1250]$. Since it was seen in Section~\ref{sec:simus} that reasonable sample size is required for good coverage an earlier start time is used here than for the estimate of the expected length of stay to ensure sufficient numbers of patients are under observation. Since we expect $P_{21}(500,t) < 0.5$, we chose the confidence bands using the transformation $\phi_2(p) = \log(-\log(1-p)).$

Figure \ref{figure_nonmarkovband} shows the confidence bands for the placebo and prednisone groups using Hall-Wellner and equal precision intervals using the non-Markov estimator and using the equal precision intervals for the Aalen-Johansen estimator. 

\begin{figure}
\begin{center}
\includegraphics[scale=0.5]{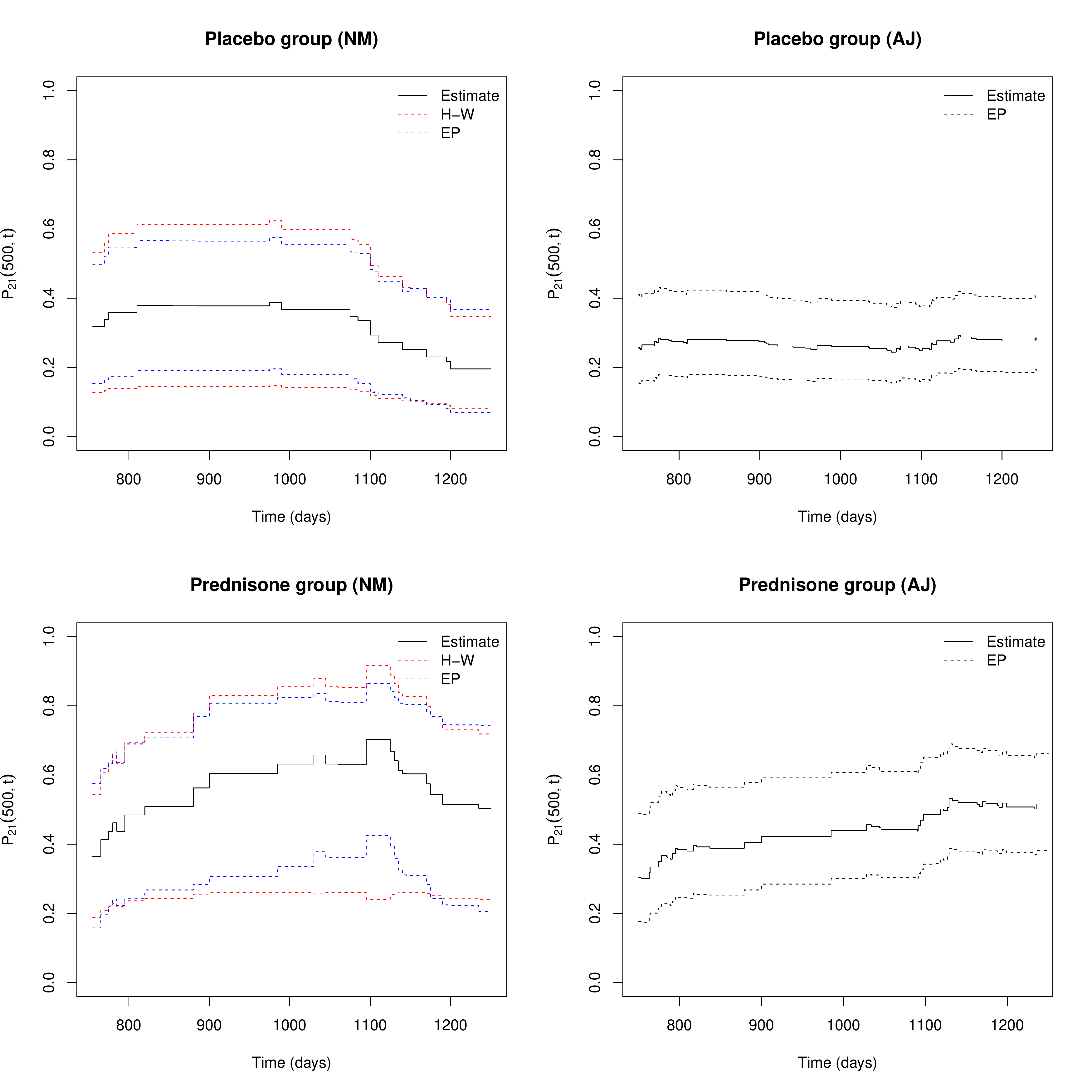}
\caption{Simultaneous confidence bands for the placebo and prednisone groups for the liver cirrhosis dataset using the non-Markov estimator (left panels) and Aalen-Johansen estimator (right panels).}
\label{figure_nonmarkovband}
\end{center}
\end{figure}

}

\section{Discussion}
\label{sec:disc}

To the best of our knowledge, our paper is the first to prove weak convergence properties of a general transition probability estimator in independently right-censored non-Markov multi-state models
while also providing explicit formulas of the asymptotic covariance functions.
Similar proofs have been given for the classical bootstrap, which finally allows the utilization of the \cite{titman15} estimator in time-simultaneous inference procedures.
Unfortunately, when focusing on inference on the expected length of stay for single time points, 
both bootstrap methods applied in Section~\ref{sec:simus} yielded a greater deviance of the nominal coverage level 
than the simple studentization method based on the asymptotic quantiles of the standard normal distribution.
Even though the simulated small sample coverage probabilities were more satisfying,
other resampling procedures might improve the performance in both, small and large samples.
For example, the martingale residual multiplier methods of \cite{lin97}, \cite{beyersmann12b}, or \cite{dobler15b} are the method of choice in incomplete competing risks data 
and they may be adapted for the present context in a future paper.

{\blue 
On the other hand, the time-simultaneous confidence bands derived in this paper cannot rely on a purely asymptotic quantile finding approach.
Instead, they truly require a resampling method such as the bootstrap method developed here.
Furthermore, the simulation results in Section~\ref{sec:bands} and the application to the liver cirrhosis data set in Section~\ref{sec:example} demonstrate their practical usefulness:
not only did they reveal reliable coverage probabilities for relatively small sample sizes already.
They were also only slightly wider than the confidence bands based on the Aalen-Johansen estimator for which there is no guarentee of applicability.
}

  Another field of possible future interest is the extension of the present methodology to quality-adjusted sums of expected lengths of stay in a specific sub-set of states;
cf. \cite{williams85} who proposed such quality-adjusted life years for operation considerations.
For instance, the expected length of stay in a healthy state could be weighted with a factor $p \in (0.5, 1)$
whereas the length of stay in the illness state obtains the weight $1 - p$.
The specific choice of $p$ may be obtained from additional information on, e.g., the pain scores of patients in the illness state and may be chosen individually by the patients themselves.
This way, a comparison of treatments could be achieved that more realistically accounts for the circumstances of a specific disease.
However, one should bear in mind the controversial debate on quality-adjusted life years as articulated by e.g. \cite{harris87};
see also the discussion of the paper by \cite{cox92}.
Therein, G. W. Torrance appreciates the usefulness of quality-adjusted life years for resource allocation if combined with other measures and advices:
{\color{black}
He states that quality-adjusted life years ``were never intended for clinical decision-making -
 they were developed for use in resource allocation.'' 
 They should also be reported together with other important measures for reaching a decision.
 Finally, one should add ``a thorough sensitivity analysis and thoughtful discussion, including caveats''
 to obtain an informative quality-adjusted life years analysis.
}

Finally, reconsider the data analysis in Section~\ref{sec:example} where the Aalen--Johansen estimator yielded a difference of expected lengths of stay
which drastically deviates from the same quantity based on the \cite{titman15} estimator.
This might be considered as a strong evidence against the Markov assumption.
If, on the other hand, a formal test of Markovianity does not find such evidence, 
the resulting estimated lengths of stay might be comparable as well.
Even more can be gained {\blue by means of such tests}: as apparent from the mentioned liver cirrhosis data example, 
the confidence intervals {\blue and simultaneous confidence bands} based on the Aalen--Johansen estimator are much narrower,
resulting from the efficiency of this estimator in comparison to non-Markov transition probability estimates.
Equivalently, test procedures based on the classical estimator are more powerful in detecting deviances from null hypotheses.
Therefore, formal tests for the applicability of the classical Aalen--Johansen estimator shall be developed in future articles.

\section*{Acknowledgements}

The authors would like to thank Markus Pauly for helpful discussions.
Dennis Dobler was supported by a DFG (German Research Foundation) grant.

\bibliographystyle{plainnat}
\bibliography{literatur}

\begin{thebibliography}{31}
\providecommand{\natexlab}[1]{#1}
\providecommand{\url}[1]{\texttt{#1}}
\expandafter\ifx\csname urlstyle\endcsname\relax
  \providecommand{\doi}[1]{doi: #1}\else
  \providecommand{\doi}{doi: \begingroup \urlstyle{rm}\Url}\fi

\bibitem[Aalen(1978)]{aalen78}
O.~O. Aalen.
\newblock {Nonparametric Inference for a Family of Counting Processes}.
\newblock \emph{The Annals of Statistics}, 6\penalty0 (4):\penalty0 701--726,
  1978.

\bibitem[Akritas(1986)]{akritas86}
M.~G. Akritas.
\newblock {B}ootstrapping the {K}aplan--{M}eier {E}stimator.
\newblock \emph{Journal of the American Statistical Association}, 81\penalty0
  (396):\penalty0 1032--1038, 1986.

\bibitem[Allignol et~al.(2014)Allignol, Beyersmann, Gerds, and
  Latouche]{allignol14}
A.~Allignol, J.~Beyersmann, T.~Gerds, and A.~Latouche.
\newblock A competing risks approach for nonparametric estimation of transition
  probabilities in a non-{M}arkov illness-death model.
\newblock \emph{Lifetime Data Analysis}, 20\penalty0 (4):\penalty0 495--513,
  2014.

\bibitem[Andersen et~al.(1993)Andersen, Borgan, Gill, and Keiding]{abgk93}
P.~K. Andersen, {\O}.~Borgan, R.~D. Gill, and N.~Keiding.
\newblock \emph{Statistical Models Based on Counting Processes}.
\newblock Springer, New York, 1993.

\bibitem[Beyersmann et~al.(2013)Beyersmann, Pauly, and Termini]{beyersmann12b}
J.~Beyersmann, M.~Pauly, and S.~Di Termini.
\newblock Weak {C}onvergence of the {W}ild {B}ootstrap for the
  {A}alen--{J}ohansen {E}stimator of the {C}umulative {I}ncidence {F}unction of
  a {C}ompeting {R}isk.
\newblock \emph{Scandinavian Journal of Statistics}, 40\penalty0 (3):\penalty0
  387--402, 2013.

\bibitem[Billingsley(1999)]{billingsley99}
P.~Billingsley.
\newblock \emph{{Convergence of Probability Measures}}.
\newblock John Wiley \& Sons, New York, second edition, 1999.

\bibitem[Bluhmki et~al.(2018)Bluhmki, Schmoor, Dobler, Pauly, Finke,
  Schumacher, and Beyersmann]{bluhmki18}
T.~Bluhmki, C.~Schmoor, D.~Dobler, M.~Pauly, J.~Finke, M.~Schumacher, and
  J.~Beyersmann.
\newblock {A wild bootstrap approach for the Aalen--Johansen estimator}.
\newblock \emph{Biometrics, early view}, 74\penalty0 (3):\penalty0 977--985,
  2018.

\bibitem[Cox et~al.(1992)Cox, Fitzpatrick, Fletcher, Gore, Spiegelhalter, and
  Jones]{cox92}
D.~R. Cox, R.~Fitzpatrick, A.~E. Fletcher, S.~M. Gore, D.~J. Spiegelhalter, and
  D.~R. Jones.
\newblock {Quality-of-life Assessment: Can We Keep It Simple?}
\newblock \emph{Journal of the Royal Statistical Society. Series A (Statistics
  in Society)}, 155\penalty0 (3):\penalty0 353--393, 1992.

\bibitem[Datta and Satten(2001)]{datta01}
S.~Datta and G.~A. Satten.
\newblock {V}alidity of the {A}alen--{J}ohansen estimators of stage occupation
  probabilities and {N}elson--{A}alen estimators of integrated transition
  hazards for non-{M}arkov models.
\newblock \emph{Statistics \& Probability Letters}, 55\penalty0 (4):\penalty0
  403--411, 2001.

\bibitem[Datta and Satten(2002)]{datta02}
S.~Datta and G.~A. Satten.
\newblock {E}stimation of {I}ntegrated {T}ransition {H}azards and {S}tage
  {O}ccupation {P}robabilities for {N}on-{M}arkov {S}ystems {U}nder {D}ependent
  {C}ensoring.
\newblock \emph{Biometrics}, 58\penalty0 (4):\penalty0 792--802, 2002.

\bibitem[de~U{\~n}a-{\'A}lvarez and Meira-Machado(2015)]{deUna-alvarez15}
J.~de~U{\~n}a-{\'A}lvarez and L.~Meira-Machado.
\newblock {Nonparametric Estimation of Transition Probabilities in the
  Non-Markov Illness-Death Model: A Comparative Study}.
\newblock \emph{Biometrics}, 71\penalty0 (2):\penalty0 364--375, 2015.

\bibitem[Dobler(2016)]{dobler16phd}
D.~Dobler.
\newblock \emph{{Nonparametric inference procedures for multi-state Markovian
  models with applications to incomplete life science data}}.
\newblock PhD thesis, Universit\"at Ulm, Deutschland, 2016.

\bibitem[Dobler and Pauly(2014)]{dobler14}
D.~Dobler and M.~Pauly.
\newblock {B}ootstrapping {A}alen-{J}ohansen processes for competing risks:
  {H}andicaps, solutions, and limitations.
\newblock \emph{Electronic Journal of Statistics}, 8\penalty0 (2):\penalty0
  2779--2803, 2014.

\bibitem[Dobler et~al.(2017)Dobler, Beyersmann, and Pauly]{dobler15b}
D.~Dobler, J.~Beyersmann, and M.~Pauly.
\newblock {N}on-strange weird resampling for complex survival data.
\newblock \emph{Biometrika}, 104\penalty0 (3):\penalty0 699--711, 2017.

\bibitem[Gill(1989)]{gill89}
R.~D. Gill.
\newblock {N}on- and {S}emi-{P}arametric {M}aximum {L}ikelihood {E}stimators
  and the von {M}ises {M}ethod ({P}art 1).
\newblock \emph{Scandinavian Journal of Statistics}, 16\penalty0 (2):\penalty0
  97--128, 1989.

\bibitem[Glidden(2002)]{glidden02}
D.~V. Glidden.
\newblock {R}obust {I}nference for {E}vent {P}robabilities with {N}on-{M}arkov
  {E}vent {D}ata.
\newblock \emph{Biometrics}, 58\penalty0 (2):\penalty0 361--368, 2002.

\bibitem[Grand and Putter(2016)]{grand16}
M.~Klinten Grand and H.~Putter.
\newblock Regression models for expected length of stay.
\newblock \emph{Statistics in Medicine}, 35\penalty0 (7):\penalty0 1178--1192,
  2016.

\bibitem[Harris(1987)]{harris87}
J.~Harris.
\newblock {QALYfying the value of life}.
\newblock \emph{Journal of medical ethics}, 13\penalty0 (3):\penalty0 117--123,
  1987.

\bibitem[Jacod and Shiryaev(2003)]{jacod03}
J.~Jacod and A.~N Shiryaev.
\newblock \emph{Limit Theorems for Stochastic Processes}.
\newblock Springer, Berlin, second edition, 2003.

\bibitem[Lin(1997)]{lin97}
D.~Y. Lin.
\newblock Non-parametric inference for cumulative incidence functions in
  competing risks studies.
\newblock \emph{Statistics in Medicine}, 16\penalty0 (8):\penalty0 901--910,
  1997.

\bibitem[Meilijson(1972)]{meilijson72}
I.~Meilijson.
\newblock {Limiting Properties of the Mean Residual Lifetime Function}.
\newblock \emph{The Annals of Mathematical Statistics}, 43\penalty0
  (1):\penalty0 354--357, 1972.

\bibitem[Meira-Machado et~al.(2006)Meira-Machado, de~U{\~n}a-{\'A}lvarez, and
  Cadarso-Su{\'a}rez]{meira-machado06}
L.~Meira-Machado, J.~de~U{\~n}a-{\'A}lvarez, and C.~Cadarso-Su{\'a}rez.
\newblock {Nonparametric estimation of transition probabilities in a non-Markov
  illness--death model}.
\newblock \emph{Lifetime Data Analysis}, 12\penalty0 (3):\penalty0 325--344,
  2006.

\bibitem[Pepe(1991)]{pepe1991}
M.~S. Pepe.
\newblock {Inference for Events with Dependent Risks in Multiple Endpoint
  Studies}.
\newblock \emph{Journal of the American Statistical Association}, 86\penalty0
  (415):\penalty0 770--778, 1991.

\bibitem[Pepe et~al.(1991)Pepe, Longton, and Thornquist]{pepe1991b}
M.~S. Pepe, G.~Longton, and M.~Thornquist.
\newblock {A qualifier Q for the survival function to describe the prevalence
  of a transient condition}.
\newblock \emph{Statistics in Medicine}, 10\penalty0 (3):\penalty0 413--421,
  1991.

\bibitem[Pollard(1984)]{pollard84}
D.~Pollard.
\newblock \emph{{Convergence of Stochastic Processes}}.
\newblock Springer, New York, 1984.

\bibitem[Putter and Spitoni(2018)]{putter16n}
H.~Putter and C.~Spitoni.
\newblock {Non-parametric estimation of transition probabilities in non-Markov
  multi-state models: The landmark Aalen--Johansen estimator}.
\newblock \emph{Statistical Methods in Medical Research}, 27\penalty0
  (7):\penalty0 2081--2092, 2018.

\bibitem[Rodr{\'\i}guez-Girondo and U{\~n}a-{\'A}lvarez(2012)]{rodriguez12}
M.~Rodr{\'\i}guez-Girondo and J.~U{\~n}a-{\'A}lvarez.
\newblock {A nonparametric test for Markovianity in the illness-death model}.
\newblock \emph{Statistics in Medicine}, 31\penalty0 (30):\penalty0 4416--4427,
  2012.

\bibitem[Titman(2015)]{titman15}
A.~C. Titman.
\newblock {Transition Probability Estimates for Non-Markov Multi-State Models}.
\newblock \emph{Biometrics}, 71\penalty0 (4):\penalty0 1034--1041, 2015.

\bibitem[van~der Vaart and Wellner(1996)]{vaart96}
A.~W. van~der Vaart and J.~A. Wellner.
\newblock \emph{Weak Convergence and Empirical Processes}.
\newblock Springer, New York, 1996.

\bibitem[Wangler et~al.(2006)Wangler, Beyersmann, and Schumacher]{wangler06}
M.~Wangler, J.~Beyersmann, and M.~Schumacher.
\newblock {changeLOS: An R-package for change in length of hospital stay based
  on the Aalen--Johansen estimator}.
\newblock \emph{R News}, 6\penalty0 (2):\penalty0 31--35, 2006.

\bibitem[Williams(1985)]{williams85}
A.~Williams.
\newblock Economics of coronary artery bypass grafting.
\newblock \emph{British Medical Journal (Clinical Research Edition)},
  291\penalty0 (6491):\penalty0 326--329, 1985.

\end{thebibliography}

\appendix

This appendix includes the proofs of Lemma~\ref{lemma:main} and the main Theorems~\ref{thm:main} to~\ref{thm:elos} as well as a detailed derivation of the limit covariance functions.
{\blue Furthermore, it is proven that the model in the simulation Section~\ref{sec:simus} satisfies the main conditions~\eqref{eq:main_cond_1}--\eqref{eq:main_cond_4}.}

\section{Asymptotic covariance functions and comments on their estimation}
{\blue We begin by stating all covariance functions which are involved in the limit theorems.
In particular, introducing the following covariance functions will prove to be useful:
\begin{align*}
 \Sigma_{\J\I}(u,v) & = P( X(s) \in \I, X(u) \in \J, X(v) \in \J, Y(u \vee v) = 1 ) 
   - p_\JI(u) p_\JI(v) \\
 \Sigma_{\I}(u,v) & = p_\I(u \vee v) 
   - p_\JI(u) p_\JI(v) \\
 \Omega_{\J\I}(u,v) & = P( X(s) \in \I, X(u) \in \J, Y(u \vee v) = 1 ) 
   - p_\JI(u) p_\I(v).
\end{align*}
With these abbreviations, we conveniently introduce the asymptotic covariance functions which appear in Lemma~1 of the manuscript.
Their derivation in detail will be done in the next section.
}
   \begin{align*}
    & cov(G_\ji(u), G_\ji(v)) \\
    & = \frac{\Sigma_\JI(u,v) + \Sigma_\I(u,v) p_\ji(u) p_\ji(v) - \Omega_\JI(u,v) p_\ji(v) - \Omega_\JI(v,u) p_\ji(u) }{p_\I(u) p_\I(v)}  \\
    & cov(L_0(r), G_\ji(u)) = \frac{1\{u < r\}}{p_\I(u)} ( P_{\I\J0}(s,u,r) - F_0(r) p_{\ji}(u) ), \\
    & cov(L_1(t), G_\ji(u)) \\
    & = \frac{1\{u < t\}}{p_\I(u)} ( P_{\I\J1}(s,u,t) - F_1(u) - p_{\ji}(u) ( F_1(t) - F_1(u) )), \\
    & cov(L_0(r), L_1(t)) = - F_0(r) \int_s^{r \wedge t} {\blue \frac{\d F_1(v)}{p_\I(v)}}.
   \end{align*}
   The asymptotic variance-covariance functions of the Kaplan--Meier and the Aalen--Johansen estimator
   for different time points are well-known and therefore {\blue not listed}.
   The derivation of $cov(L_0(r), G_\ji(u))$ and $cov(L_1(t), G_\ji(u))$ will be made explicit in Section~\ref{sec:cov} below.
   {\blue In the display above, we used} $P_{\I\J k}(s,u,t) = P(Z(t) = k, X(u) \in \J \mid X(s) \in \I ), k=1,2$.
Note that {\blue these} functions may be decomposed into
\begin{align*}
 P( Z(t) = k \mid X(u) \in \J, X(s) \in \I ) P_{\I\J}(s,u)
\end{align*}
of which the first probability is again estimable by Kaplan--Meier or Aalen--Johan\-sen-type estimators
based on those individuals that are observed to satisfy $X(u) \in {\blue \J}$ and $X(s) \in \I$.
{\blue Hence, consistent plug-in estimates of the above covariance function are found easily.}
However, the number of individuals satisfying {\blue that $X(u) \in {\blue \J}$ and that $X(s) \in \I$ is observed} might be quite small {\blue in practical applications. This could lead to a large variance of the estimators of} $P_{\I\J k}(s,u,t)$.

\section{Proof of Lemma 1}
\label{sec:proofL1}

The large sample properties concerning the conditional Kaplan--Meier estimator $\wh F_0$ 
and the conditional Aalen--Johansen estimator $\wh F_1$ follow from standard theory.
Their variance-covariance functions are obtained in the same way.
The asymptotic covariance functions between $\wh p_\ji$ and $\wh F_0$ as well as $\wh F_1$
are derived in Appendix~\ref{sec:cov} below.
For technical reasons, we consider weak convergence on the c\`adl\`ag space $D[s,\tau]$
equipped with the supremum norm instead of the usual Skorohod metric.
In case of continuous limit sample paths, weak convergence on the Skorohod and the supremum-normed c\`adl\`ag space are equivalent; cf. \cite{gill89}, p. 110 in \cite{abgk93}, or p. 137 in \cite{pollard84}.
The advantage is that the functional delta-method is available on the supremum-normed space.
Additionally, the general Aalen--Johansen estimator is known to have continuous limit distribution sample paths; cf. Theorem~IV.4.2 in \cite{abgk93}.
However, this may not be true for the weak limit of $\sqrt{n}(\wh p_\ji - p_\ji)$ as $n \rightarrow \infty$
because the censoring distribution may have discrete components.
This problem is solved by an application of a technique as used by \cite{akritas86}:
Each discrete censoring components is distributed uniformly on adjacent, inserted small time intervals 
during which no state transition can occur.
Instead, all future state transitions are shifted by the preceding interval lengths.
Write $ W_{\J\I} =  \sqrt{n} ( \wh p_{\J\I} - p_{\J\I})$ 
  and $ W_\I = \sqrt{n} ( \wh p_{\I} - p_\I)$,
  where $p_{\JI}(t) = P(X(s) \in \I, X(t) \in \J, Y(t) = 1)$ 
  and $p_{\I}(t) = P(X(s) \in \I, Y(t) = 1)$.
The resulting weak limit sample paths of the time transformed $\sqrt{n}(\wh p_\I - p_\I)$ and $\sqrt{n}(\wh p_\JI - p_\JI)$ will be continuous and the functional delta-method is applicable to obtain the corresponding result for $\wh p_\ji$.
Note that the values of none of these processes were modified outside the just inserted small time intervals.
Finally, the projection functional, which projects a time-transformed process to the original time line, is continuous.
Hence, the continuous mapping theorem eventually yields the desired central limit theorem for the process
$\sqrt{n}(\wh p_\ji - p_\ji)$ as $n \rightarrow \infty$.

It remains to give a weak convergence proof of the finite-dimensional margins of the process 
$\sqrt{n} (\wh p_\ji - p_\ji)$ and to verify its tightness under the temporary assumption of continuously distributed censoring times.
Combining both properties, weak convergence on $D[s,\tau]$ follows, as indicated above.

\subsection*{Finite-Dimensional Marginal Distributions}

Weak convergence of the finite-dimensional margins is shown by applying the functional delta-method
after utilizations of the central limit theorem.
  The weak convergences of all finite-dimensional marginal distributions are easily 
  obtained by the multivariate central limit theorem:
  For $m \in \N$ points of time, $s \leq t_1 < \dots < t_m \leq \tau$, consider
  $(W_{\JI}(t_1), \dots, W_{\JI}(t_m))' \rightarrow (\tilde G_{\JI}(t_1), \dots, \tilde G_{\JI}(t_m))' \sim N(0,\Sigma_{\JI})$ in distribution, as $n \rightarrow \infty$,
  with $(j, k)$th covariance entry
  \begin{align*}
     P( X(s) \in \I, X(t_j) \in \J, X(t_k) \in \J, Y(t_j \vee t_k) = 1 ) 
   - p_\JI(t_j) p_\JI(t_k).
  \end{align*}
  Similarly, 
  $(W_{\I}(t_1), \dots, W_{\I}(t_m))' \rightarrow (\tilde G_{\I}(t_1), \dots, \tilde G_{\I}(t_m))' \sim N(0,\Sigma_{\I})$ in distribution, as $n \rightarrow \infty$,
  with $(j, k)$th covariance entry 
  \begin{align*}
   & cov(W_\I(t_j), W_\I(t_k)) = p_\I(t_j \vee t_k) - p_\I(t_j) p_\I(t_k).
  \end{align*}
  The convergence of both vectors even holds jointly,
  implying that the $(j,k)$th entry of the covariance matrix $\Omega_\JI$ of both limit vectors is given by
  \begin{align*}
   \Omega^{[j,k]}_{\JI} =& cov(\tilde G_{\JI}(t_j), \tilde G_\I(t_k)) 
   = cov(W_{\JI}(t_j), W_\I(t_k)) \\
   & = P( X(s) \in \I, X(t_j) \in \J, Y(t_j \vee t_k) = 1 ) - p_\JI(t_j) p_\I(t_k).
  \end{align*}
  Denote by $D_{> 0}[s,\tau]$ the subset of c\`adl\`ag functions with are positive and bounded away from zero, which is equipped with the supremum norm.
  Note that $\phi: D[0,\tau] \times D_{> 0}[s,\tau] \rightarrow D[s,\tau],  (f,g) \mapsto \frac{f}g$ defines a Hadamard-differentiable map
  with Hadamard-derivative $\phi'_{(f,g)} (h_1,h_2) = \frac{h_1}{g} - h_2 \frac{f}{g^2}$, with $h_1, h_2 \in D[s,\tau]$.
  This of course also implies Fr\'echet-differentiability on all finite-dimensional projections.
  Therefore, the multivariate delta-method yields the convergence in distribution
  \begin{align*}
   \sqrt{n} & (\wh p_{\J|\I}(t_1) - p_{\J|\I}(t_1), \dots, \wh p_{\J|\I}(t_m) - p_{\ji}(t_m))' 
     \rightarrow (G_{\ji}(t_1), \dots, G_{\ji}(t_m))' \\
    & = \Big(\frac{\tilde G_{\JI}(t_1)}{p_{\I}(t_1)} - \tilde G_{\I}(t_1) \frac{p_{\ji}(t_1)}{p^2_{\I}(t_1)}, \dots, \frac{\tilde G_{\JI}(t_m)}{p_{\I}(t_m)} - \tilde G_{\I}(t_m) \frac{p_{\JI}(t_m)}{p^2_{\I}(t_m)} \Big)',
  \end{align*}
  as $n \rightarrow \infty$.
  The asymptotic covariance matrix $\Gamma_{\ji}$ has $(j,k)$th entry
  \begin{align*}
   & \frac{\Sigma_{\JI}^{[j,k]}}{p_\I(t_j) p_\I(t_k)} + {\Sigma_{\I}^{[j,k]}} \frac{p_{\JI}(t_j) p_{\JI}(t_k)}{p_\I^2(t_j) p_\I^2(t_k)}
   - \Omega^{[j,k]}_{\JI} \frac{p_{\JI}(t_k)}{p_\I(t_j) p_\I^2(t_k)} 
   - \Omega^{[k,j]}_{\JI} \frac{p_{\JI}(t_j)}{p_\I(t_k) p_\I^2(t_j)} \\
   & = \frac{1}{p_\I(t_j) p_\I(t_k)} \Big( \Sigma_{\JI}^{[j,k]} 
    + \Sigma_{\I}^{[j,k]} p_{\ji}(t_j) p_{\ji}(t_k)
    - \Omega^{[j,k]}_{\JI} p_{\ji}(t_k) - \Omega^{[k,j]}_{\JI} p_{\ji}(t_j) \Big) .
  \end{align*}

\subsection*{Tightness}

We will show tightness on the Skorohod space $D[s,\tau]$ of c\`adl\`ag functions 
by showing that a modulus of continuity becomes small in probability; see Theorem~13.5 in \cite{billingsley99}.
Remember that this property will imply weak convergence on the Skorohod space $D[s,\tau]$
which, due to the temporary continuity of the limit processes, is equivalent to weak convergence on the supremum-normed space.

  Now, for the modulus of continuity,
  we consider for any small $\delta >0$ the expectation
  \begin{align}
  \label{eq:mod_W}
   & \sup_{\substack{s \leq t_1 \leq t \leq t_2 \leq \tau \\ t_2 - t_1 \leq \delta}}
    E \Big( |W_\I(t) - W_\I(t_1)|^2 |W_\I(t_2) - W_\I(t)|^2  \Big).
  \end{align}
  We start by rewriting $W_\I(t) - W_\I(t_1)$ as
  \begin{align*}
  & \frac{1}{\sqrt n} \sum_{\ell=1}^n ( 1 \{ X_\ell(s) \in \I, Y_\ell (t_1) = 1, Y_\ell (t) = 0 \}
   - P( X(s) \in \I,  Y (t_1) = 1, Y (t) = 0 ) ) \\
  & = \frac{1}{\sqrt n} \sum_{\ell=1}^n (Z_{\I10; \ell}(s,t_1,t) - r_{\I10}(s,t_1,t)).
  \end{align*}
  Note that, in the product appearing in \eqref{eq:mod_W}, the products of the above indicator functions with identical indices vanish for unequal time pairs $(t_1,t)$ and $(t, t_2)$.
  Additionally, due to the independence of all individuals, all terms with a single separate index $\ell$ vanish too.
  Hence, by an application of the Cauchy-Schwarz inequality, the expectation in \eqref{eq:mod_W} is bounded above by
  \begin{align*}
   3 & E[ (Z_{\I10; 1}(s,t_1,t) - r_{\I10}(s,t_1,t))^2] 
   E[(Z_{\I10; 1}(s,t,t_2) - r_{\I10}(s,t,t_2))^2] \\
   & + \frac1n E[(Z_{\I10; 1}(s,t_1,t) - r_{\I10}(s,t_1,t))^2 (Z_{\I10; 1}(s,t,t_2) - r_{\I10}(s,t,t_2))^2 ] \\
    = & 3 r_{\I10}(s,t_1,t) ( 1 - r_{\I10}(s,t_1,t) ) r_{\I10}(s,t,t_2) (1- r_{\I10}(s,t,t_2) ) \\
   & + \frac1n [ r^2_{\I10}(s,t_1,t) r_{\I10}(s,t,t_2) ( 1-  2 r_{\I10}(s,t,t_2) ) \\
   & + r^2_{\I10}(s,t,t_2) r_{\I10}(s,t_1,t) ( 1 - 2 r_{\I10}(s,t_1,t) ) 
    +  r^2_{\I10}(s,t_1,t) r^2_{\I10}(s,t,t_2) ] \\
    \leq & 4 r_{\I10}(s,t_1,t) r_{\I10}(s,t,t_2) 
  \end{align*}
  for all $n \geq 3$. 
  Let $S_C(t) = \pr( C > t)$ denote the survival function of the censoring times
  and let $r_1(t) = P(Y(t) = 1) \leq S_C(t)$ denote the (unconditional) probability to be still under risk at time $t$.
  By the temporarily continuous censoring distribution, 
  $r_1$ is a continuous and non-decreasing function and thus it is suitable for an application of Theorem~13.5 in \cite{billingsley99}.
  Then the terms in the above display are not greater than
  $$ 4 (r_1(t_1) - r_1(t)) (r_1(t) - r_1(t_2)) \leq 4 (r_1(t_1) - r_1(t_2))^2. $$
  Hence, tightness in the Skorohod space $D[s,\tau]$ of the process $\sqrt{n}(\wh p_\I - p_\I)$ affiliated with the denominator of $\wh p_{\ji}$ follows.

  It remains to prove tightness of the numerator process.
  At first, we rewrite for $t_1 < t < t_2$
  \begin{align*}
   W_{\JI}(t) & - W_{\JI}(t_1) 
   =  \frac{1}{\sqrt{n}} \sum_{\ell=1}^n (\zeta_{\I01;\ell}(s,t_1,t) 
    - \zeta_{\I10;\ell}(s,t_1,t)
    - \zeta_{\I11;\ell}(s,t_1,t) \\
   & - q_{\I01}(s,t_1,t) + q_{\I10}(s,t_1,t) + q_{\I11}(s,t_1,t)) \\
   = & \frac{1}{\sqrt{n}} \sum_{\ell = 1}^n
   ( 1 \{ X_\ell(s) \in \I, X_\ell(t_1) \notin J, X_\ell(t) \in \J, Y_\ell (t) = 1 \} \\
   & - 1 \{ X_\ell(s) \in \I, X_\ell(t_1) \in \J, X_\ell(t) \notin \J, Y_\ell (t) = 1 \} \\
   & - 1 \{ X_\ell(s) \in \I, X_\ell(t_1) \in \J, Y_\ell (t_1) = 1, X_\ell(t) \in \J, Y_\ell (t) = 0 \} \\
   & - E(\zeta_{\I01;1}(s,t_1,t)) + E(\zeta_{\I10;1}(s,t_1,t)) + E(\zeta_{\I11;1}(s,t_1,t)) ).  
  \end{align*}
  For later calculations, we derive products of the $\zeta$-indicators for different points of time:
  \begin{align*}
   & \zeta_{\I01;\ell}(s,t_1,t) \zeta_{\I01;\ell}(s,t,t_2) = 0, \\
   & \zeta_{\I10;\ell}(s,t_1,t) \zeta_{\I10;\ell}(s,t,t_2) = 0, \\
   & \zeta_{\I11;\ell}(s,t_1,t) \zeta_{\I11;\ell}(s,t,t_2) = 0, \\
   & \zeta_{\I01;\ell}(s,t_1,t) \zeta_{\I10;\ell}(s,t,t_2) \\
   & = 1 \{ X_\ell(s) \in \I, X_\ell(t_1) \notin \J, X_\ell(t) \in \J, X_\ell(t_2) \notin \J, Y_\ell (t_2) = 1 \} \\
   & = \zeta_{\I010;\ell}(s,t_1,t,t_2), \\
   & \zeta_{\I10;\ell}(s,t_1,t) \zeta_{\I01;\ell}(s,t,t_2) \\
   & = 1 \{ X_\ell(s) \in \I, X_\ell(t_1) \in \J, X_\ell(t) \notin \J, X_\ell(t_2) \in \J, Y_\ell (t_2) = 1 \} \\
   & = \zeta_{\I101;\ell}(s,t_1,t,t_2), \\
   & \zeta_{\I01;\ell}(s,t_1,t) \zeta_{\I11;\ell}(s,t,t_2) \\
   & = 1 \{ X_\ell(s) \in \I, X_\ell(t_1) \notin \J, X_\ell(t) \in \J, X_\ell(t_2) \in \J, Y_\ell (t) = 1, Y_\ell (t_2) = 0 \} \\
   & = \zeta_{\I011;\ell}(s,t_1,t,t_2), \\
   & \zeta_{\I11;\ell}(s,t_1,t) \zeta_{\I01;\ell}(s,t,t_2) = 0, \\
   & \zeta_{\I10;\ell}(s,t_1,t) \zeta_{\I11;\ell}(s,t,t_2) = 0, \\
   & \zeta_{\I11;\ell}(s,t_1,t) \zeta_{\I10;\ell}(s,t,t_2) = 0.
  \end{align*}
  Therefore, the expectation to be bounded for verifying tightness has the following structure:
  \begin{align*}
   E & \Big( |W_{\JI}(t) - W_{\JI}(t_1)|^2 |W_{\JI}(t_2) - W_{\JI}(t)|^2  \Big) \\
   \leq & \frac1n E \Big( 
   (\zeta_{\I01;1}(s,t_1,t) - \zeta_{\I10;1}(s,t_1,t) - \zeta_{\I11;1}(s,t_1,t) \\
   & -  q_{\I01}(s,t_1,t) + q_{\I10}(s,t_1,t) + q_{\I11}(s,t_1,t))^2 \\
   & (\zeta_{\I01;1}(s,t,t_2) - \zeta_{\I10;1}(s,t,t_2) - \zeta_{\I11;1}(s,t,t_2) \\
   & - q_{\I01}(s,t,t_2) + q_{\I10}(s,t,t_2) + q_{\I11}(s,t,t_2))^2 \Big) \\
   & + 3 E \Big( 
   (\zeta_{\I01;1}(s,t_1,t) - \zeta_{\I10;1}(s,t_1,t) - \zeta_{\I11;1}(s,t_1,t) \\
   & -  q_{\I01}(s,t_1,t) + q_{\I10}(s,t_1,t) + q_{\I11}(s,t_1,t))^2 \Big) \\
   & E \Big( (\zeta_{\I01;1}(s,t,t_2) - \zeta_{\I10;1}(s,t,t_2) - \zeta_{\I11;1}(s,t,t_2) \\
   & - q_{\I01}(s,t,t_2) + q_{\I10}(s,t,t_2) + q_{\I11}(s,t,t_2))^2 \Big),
  \end{align*}
  where the inequality follows from an application of the Cauchy-Schwarz inequality.
  The terms in the previous display are not greater than
  \begin{align*}
   & \frac9n E \Big( [
   (\zeta_{\I01;1}(s,t_1,t) -  q_{\I01}(s,t_1,t))^2 
    + (\zeta_{\I10;1}(s,t_1,t) - q_{\I10}(s,t_1,t))^2 \\
   & + (\zeta_{\I11;1}(s,t_1,t) - q_{\I11}(s,t_1,t))^2 ]
    [ (\zeta_{\I01;1}(s,t,t_2) - q_{\I01}(s,t,t_2))^2 \\
   & + (\zeta_{\I10;1}(s,t,t_2) - q_{\I10}(s,t,t_2))^2 
    + (\zeta_{\I11;1}(s,t,t_2) - q_{\I11}(s,t,t_2))^2 ] \Big) \\ 
   & + 27 E \Big( 
   (\zeta_{\I01;1}(s,t_1,t) -  q_{\I01}(s,t_1,t))^2 
    + (\zeta_{\I10;1}(s,t_1,t) - q_{\I10}(s,t_1,t))^2 \\
   & + (\zeta_{\I11;1}(s,t_1,t) - q_{\I11}(s,t_1,t))^2  \Big)
   E \Big(  (\zeta_{\I01;1}(s,t,t_2) - q_{\I01}(s,t,t_2))^2 \\
   & + (\zeta_{\I10;1}(s,t,t_2) - q_{\I10}(s,t,t_2))^2 
    + (\zeta_{\I11;1}(s,t,t_2) - q_{\I11}(s,t,t_2))^2  \Big) .
  \end{align*}
  It remains to calculate and bound the single expectations.
  \begin{align*}
   E \Big(   (\zeta_{\I01;1}(s,t_1,t) -  q_{\I01}(s,t_1,t))^2 \Big)
   = q_{\I01}(s,t_1,t) (1 - q_{\I01}(s,t_1,t)) \leq q_{\I01}(s,t_1,t)
  \end{align*}
  holds by the same reasons as above and a similar representation holds for the other second moments.
  By conditions~(1) to (4) in the manuscript, 
  there is a global constant $K > 0$, 
  that is independent of $t_1,t,t_2,n$,
  such that $ q_{\I01}(s,t_1,t) \leq ({\blue \lambda_{\bar\J\J}(t)} + K)(t - t_1)$.
  Similarly, $ q_{\I10}(s,t_1,t) \leq ({\blue \lambda_{\J\bar\J}(t)} + K)(t - t_1)$
  as well as $ q_{\I11}(s,t_1,t) \leq r_1(t_1) - r_1(t)$.

  The first expectation (involving fourth moments) contains the following terms:
  \begin{align*}
   e_1 = & E \Big(  (\zeta_{\I01;1}(s,t_1,t) -  q_{\I01}(s,t_1,t))^2 
   (\zeta_{\I01;1}(s,t,t_2) - q_{\I01}(s,t,t_2))^2 \Big), \\
   e_2 = & E \Big(  (\zeta_{\I10;1}(s,t_1,t) - q_{\I10}(s,t_1,t))^2 
   (\zeta_{\I10;1}(s,t,t_2) - q_{\I10}(s,t,t_2))^2 \Big), \\
   e_3 = & E \Big(  (\zeta_{\I11;1}(s,t_1,t) - q_{\I11}(s,t_1,t))^2 
   (\zeta_{\I11;1}(s,t,t_2) - q_{\I11}(s,t,t_2))^2 \Big), \\
   e_{12} = & E \Big( (\zeta_{\I01;1}(s,t_1,t) -  q_{\I01}(s,t_1,t))^2 
   (\zeta_{\I10;1}(s,t,t_2) - q_{\I10}(s,t,t_2))^2 \Big), \\
   e_{21} = & E \Big( (\zeta_{\I10;1}(s,t_1,t) -  q_{\I10}(s,t_1,t))^2 
   (\zeta_{\I01;1}(s,t,t_2) - q_{\I01}(s,t,t_2))^2 \Big), \\
   e_{13} = & E \Big( (\zeta_{\I01;1}(s,t_1,t) -  q_{\I01}(s,t_1,t))^2 
   (\zeta_{\I11;1}(s,t,t_2) - q_{\I11}(s,t,t_2))^2 \Big), \\
   e_{31} = & E \Big( (\zeta_{\I11;1}(s,t_1,t) -  q_{\I11}(s,t_1,t))^2 
   (\zeta_{\I01;1}(s,t,t_2) - q_{\I01}(s,t,t_2))^2 \Big), \\
   e_{23} = & E \Big( (\zeta_{\I10;1}(s,t_1,t) - q_{\I10}(s,t_1,t))^2
   (\zeta_{\I11;1}(s,t,t_2) - q_{\I11}(s,t,t_2))^2 \Big), \\
   e_{32} = & E \Big( (\zeta_{\I11;1}(s,t_1,t) - q_{\I11}(s,t_1,t))^2
   (\zeta_{\I10;1}(s,t,t_2) - q_{\I10}(s,t,t_2))^2 \Big).
  \end{align*}
  For the expectation, first we obtain
  \begin{align*}
   e_1 & = 0 - 2 E( \zeta_{\I01;1}(s,t_1,t) q_{\I01}(s,t_1,t) q^2_{\I01}(s,t,t_2) ) \\
    & - 2 E( \zeta_{\I01;1}(s,t,t_2) q_{\I01}(s,t,t_2) q^2_{\I01}(s,t_1,t) )
      + q^2_{\I01}(s,t_1,t) q^2_{\I01}(s,t,t_2) \\
    & = - 3 q^2_{\I01}(s,t_1,t) q^2_{\I01}(s,t,t_2).
  \end{align*}
  Similarly, we have $e_2, e_3, e_{31}, e_{23}, e_{32} \leq 0 $.
  It remains to consider $e_{12}, e_{13}$, and $e_{21}$.
  \begin{align*}
   e_{21} & = E( \zeta_{\I10;1}(s,t_1,t) (\zeta_{\I10;1}(s,t_1,t) - 2  q_{\I10}(s,t_1,t)) \\
     & \times \zeta_{\I01;1}(s,t,t_2) ( \zeta_{\I01;1}(s,t,t_2) - 2 q_{\I01}(s,t,t_2)) )
    - 3 q^2_{\I10}(s,t_1,t) q^2_{\I01}(s,t,t_2) \\
    & = E(\zeta_{\I101;1}(s,t_1,t,t_2)) (1 - 2 q_{\I01}(s,t,t_2) - 2 q_{\I10}(s,t_1,t) ) \\
    & - 3 q^2_{\I10}(s,t_1,t) q^2_{\I01}(s,t,t_2) \\
    & \leq P( X_1(s) \in \I, X_1(t_1) \in \J, X_1(t) \notin \J, X_1(t_2) \in \J, Y_1 (t_2) = 1 ) \\
    & \leq P( X_1(t_1) \in \J, X_1(t) \notin \J, X_1(t_2) \in \J).
  \end{align*}
  Similarly,
  $ e_{12} \leq P( X_1(t_1) \notin \J, X_1(t) \in \J, X_1(t_2) \notin \J)  $
  and $ e_{13} \leq P( X_1(t_1) \notin \J, X_1(t) \in \J, C_1 \in (t,t_2] ) 
  = P( X_1(t_1) \notin \J, X_1(t) \in \J) P(C_1 \in (t,t_2] ) $
  due to independent right-censoring.
  
  Assume that $e_{12}, e_{13}, e_{21} > 0$. Otherwise, nothing needs to be shown.
  Thus,
  \begin{align*}
   e_{12} \leq & P( X_1(t_2) \notin \J | X_1(t) \in \J, X_1(t_1) \notin \J) P( X_1(t) \in \J | X_1(t_1) \notin \J) \\
    & \quad \times P( X_1(t_1) \notin \J) \\
    \leq & P( X_1(t_2) \notin \J | X_1(t) \in \J, X_1(t_1) \notin \J) P( X_1(t) \in \J | X_1(t_1) \notin \J).
  \end{align*}
  Due to condition~(4), there is a constant $K > 0$, which is independent of $t_1,t,t_2,n$, such that
  $$ e_{12} \leq ( {\blue \lambda_{\J\bar\J\mid\bar\J}(t_2 | t )} + K ) (t_2 - t) ( {\blue\lambda_{\bar\J\J}(t)} + K ) (t - t_1) 
    \leq \tilde K^2 (t_2 - t_1)^2$$
   for a global constant $\tilde K > 0$.
   Similarly $e_{12} \leq \tilde K^2 (t_2 - t_1)^2$.
   Therefore, we conclude that 
   $$e_{13} \leq ( {\blue \lambda_{\bar\J\J}(t)} + K ) (t - t_1) ( S_C(t) - S_C(t_2)) \leq \tilde K^2 (H(t_2) - H(t_1))^2, $$
   for the non-decreasing and continuous function $H : x \mapsto x + 1 - S_C(x)$.
   
  All in all, we conclude that 
  \begin{align*}
   & E \Big( |W_{\JI}(t) - W_{\JI}(t_1)|^2 |W_{\JI}(t_2) - W_{\JI}(t)|^2  \Big) \\
   & \leq 27 \tilde K^2 (t_2 - t_1)^2 + 27 ( 2 \tilde K (t_2 - t_1) + (r_1(t_1) - r_1(t_2)))^2 \\
   & \leq ( F(t_2) - F(t_1))^2,
  \end{align*}
  where $F(x) = 9 \tilde K (x + 1 - r_1(x))$
  defines a continuous and non-decreasing function, independently of $t_1,t,t_2,n$.
  Hence, tightness of the process $W_{\JI}$ follows
  which in turn implies convergence in distribution on the Skorohod space $D[s,\tau]$.
  Finally, proceed as indicated in the beginning of this Appendix section.

\section{Proof of Theorem 1}

The proof of Theorem consists of an application of the functional delta-method applied to the $D^3[s,\tau] \rightarrow D[s,\tau]$ functional
$(f_0, f_1, p) \mapsto f_1 + f_0 p$ which is Hadamard-differentiable at $(F_0, F_1, p_\ji)$
with derivative
$$(h_0,h_1, k) \mapsto h_1 + h_0 p_\ji + k F_0.$$
For the asymptotic covariance function, note that both terms
$cov(L_0(t), G_\ji(r))$ and $cov(L_1(t), G_\ji(r))$ vanish for $r \geq t$.

\section{Proof of Theorem 2}

First, the finite-dimensional, conditional marginal distributions of the estimator $\wh P_{ij}^*$ given $\chi$ are shown to be consistent, as $n \rightarrow \infty$ for the corresponding marginal distributions of the limit process $L_1 + L_0 p_{\ji} + G_{\ji} F_0$.
  This is verified via empirical process theory:
  Let $s \leq t_1 < \dots < t_d \leq \tau$ be finitely many time points.
  The finite-dimensional marginals of the estimator $\wh P_\ji(s, \cdot)$ are Hadamard-differenti\-able functionals of the empirical process
  based on 
  \begin{align*}
   ( X_\ell(s), X_\ell(t_1), X_\ell(t_2), \dots, X_\ell(t_d), Y_\ell(s), T_\ell, Z_\ell(T_\ell), C_\ell \cdot 1\{ C_\ell < T_\ell \} ), \\
  \hfill \ell=1, \dots, n,
  \end{align*}
  where $T_\ell = \inf \{t \geq s: Y_\ell(t) \neq 1 \} \wedge \tau$.
  The empirical process is indexed by the family of functions
  $$ \mc F = \{f_0, f_1, \dots, f_d, g_0, g_{1t}, h_1, h_2, k_t : s \leq t \leq \tau \}, $$
  where the particular functions are defined as
  \begin{align*}
   f_0(x_0, x_1, \dots, x_d, y, t_0, z, c ) & = 1( x_0 \in \I), \\
   f_j(x_0, x_1, \dots, x_d, y, t_0, z, c ) & = 1( x_j \in \J), \quad j = 1, \dots, d, \\
   g_0(x_0, x_1, \dots, x_d, y, t_0, z, c ) & = 1( y = 1 ), \\
   g_{1t}(x_0, x_1, \dots, x_d, y, t_0, z, c ) & = 1( t_0 \leq t ), \\
   h_1(x_0, x_1, \dots, x_d, y, t_0, z, c ) & = 1( z = 1 ), \\
   h_2(x_0, x_1, \dots, x_d, y, t_0, z, c ) & = 1( z = 2 ), \\
   k_t(x_0, x_1, \dots, x_d, y, t_0, z, c ) & = 1( c \leq t ).
  \end{align*}
  Certainly, $\mc F$ defines a Donsker class.
  The derivation of the Kaplan--Meier estimator based on parts of this empirical process
  has been verified in Section~3.9 of \cite{vaart96}.
  This similarly holds for the Aalen--Johansen estimator under competing risks by another, final application of the Wilcoxon functional and the use of the cause-specific Nelson--Aalen estimator;
  see Section~\ref{sec:cov} below for a more detailed derivation.
  The proportion-type estimator $\wh p_\ji(t_j), j=1,\dots,d, $ is given by a simple Hadamard-differentiable functional of the empirical process indexed by the functions $f_0, f_j, g_{1t}$.
  Therefore, combining again all separate estimators, the functional delta-method is applicable to the bootstrap empirical process and the same Hadamard-differentiable functional as well (in outer probability);
  see Chapter~3.6 in \cite{vaart96} for details.
  It follows that the conditional finite-dimensional marginal distributions of $\sqrt{n} ( P^*_{\IJ}(s,\cdot) - \wh P_{\IJ}(s,\cdot))$
  converge weakly to the limit distributions of the original normalized estimator.
  
  Second, conditional tightness in probability needs to be verified.
  To this end, we utilize a variant of the tightness criterion in Theorem~13.5 in \cite{billingsley99},
  in which the involved continuous, non-decreasing function $F$ may depend on the sample size $n$ and it need not be continuous,
  but only point-wise convergent and asymptotically continuous.
  See the Appendix~A.1 in the PhD thesis by \cite{dobler16phd}
  or the comment on p. 356 in \cite{jacod03} for details.
  As in the first part of this appendix, we assume without loss of generality that the censoring distribution is continuous.
  To verify tightness, we consider the conditional version of expectation 
  that we bounded accordingly for verifying tightness of the original, normalized estimator in the first part of this appendix.
  All quantities with a superscript star are the obvious bootstrap counterparts of the quantities introduced earlier.
  We again focus on proving tightness for the bootstrap counterparts 
  $W_{\JI}^* = \sqrt{n}(p_{\JI}^* - \wh p_{\JI})$ and $W_\I^* = \sqrt{n}(p_{\I}^* - \wh p_{\I})$
  of $W_{\JI} = \sqrt{n}(\wh p_{\JI} - p_{\JI})$ and $W_\I = \sqrt{n}(\wh p_{\I} - p_{\I})$, respectively.
  Furthermore, the normalized, bootstrapped Kaplan--Meier estimator is well-known to be tight
  and the same property is easily seen for the normalized, bootstrapped Aalen--Johansen estimator under competing risks.
  Finally, we will apply the continuous mapping theorem in order to conclude asymptotic tightness of $\sqrt{n}(P_{\JI}^*(s,\cdot) - \wh P_{\JI}(s,\cdot))$.
  
  Let us abbreviate $W^*_{\JI}(t) - W^*_{\JI}(t_1) = \frac1{\sqrt n} \sum_{\ell = 1}^n (m_\ell-1) a_\ell $
  and $W^*_{\JI}(t_2) - W^*_{\JI}(t) = \frac1{\sqrt n} \sum_{\ell = 1}^n (m_\ell-1) b_\ell $,
  where $(m_1,\dots,m_n) \sim Mult(n,\frac1n, \dots, \frac1n)$  is independent of $\chi$.
  Let $s \leq t_1 \leq t \leq t_2 \leq \tau$.
  Denote by $K > 0$ a constant which may become larger along the subsequent inequalities, but which is independent of $n, t, t_1, t_2$.
  Then, by the Cauchy-Schwarz inequality,
  \begin{align*}
   E & \Big( |W^*_{ji}(t) - W^*_{ji}(t_1)|^2 |W^*_{ji}(t_2) - W^*_{ji}(t)|^2 \mid \chi \Big) \\
   \leq & \frac{K}{n^2} \Big[ \sum_{\ell=1}^n a_\ell^2 b_\ell^2 E(m_1-1)^4 
    + \sum_{\ell \neq k} a_\ell^2 b_k^2 E(m_1-1)^2(m_2-1)^2 \\
    & + \Big( \Big(\sum_{\ell=1}^n a_\ell^2\Big)^{3/2} \Big(\sum_{k=1}^n b_k^2\Big)^{1/2}
	+ \Big(\sum_{\ell=1}^n a_\ell^2\Big)^{1/2} \Big(\sum_{k=1}^n b_k^2\Big)^{3/2}\Big) E(m_1-1)^3(m_2-1) \\
    & + \Big( \sum_{\ell \neq k \neq h \neq \ell} (a_\ell^2 b_k b_h + a_\ell a_k b_h^2) 
	+ \sum_{\ell \neq k}^n a_\ell b_k \Big( \sum_{\ell=1}^n a_\ell^2 \sum_{k=1}^n b_k^2 \Big)^{1/2} \Big) \\
    & \quad \times E(m_1-1)^2 (m_2-1)(m_3-1) \\
    & + \Big( \sum_{g=1}^n \sum_{h \neq g} \sum_{\ell \neq g,h} \sum_{k \neq g,h,\ell} a_g a_h b_\ell b_k \Big) E\Big( \prod_{\ell=1}^4 (m_\ell -1) \Big) \Big] 
  \end{align*}
  As shown in e.g. \cite{dobler14}, we have $E\Big( \prod_{\ell=1}^4 (m_\ell -1) \Big) = O(n^{-2})$ and $E(m_1-1)^2 (m_2-1)(m_3-1) = O(n^{-1})$. 
  Furthermore, it is easy to show that $E(m_1-1)^4,  E(m_1-1)^2(m_2-1)^2, E(m_1-1)^3(m_2-1) < \infty$.
  Therefore, the expectation in the previous display is bounded above by
  \begin{align*}
   & \frac{K}{n^2} \Big[ \sum_{\ell=1}^n a_\ell^2 b_\ell^2
    + \sum_{\ell \neq k} a_\ell^2 b_k^2 + \Big( \Big(\sum_{\ell=1}^n a_\ell^2\Big)^{3/2} \Big(\sum_{k=1}^n b_k^2\Big)^{1/2}
	+ \Big(\sum_{\ell=1}^n a_\ell^2\Big)^{1/2} \Big(\sum_{k=1}^n b_k^2\Big)^{3/2}\Big) \\
    & + \frac1n \Big( \sum_{\ell \neq k \neq h \neq \ell} (a_\ell^2 b_k b_h + a_\ell a_k b_h^2) 
	+ \sum_{\ell \neq k}^n a_\ell b_k \Big( \sum_{\ell=1}^n a_\ell^2 \sum_{k=1}^n b_k^2 \Big)^{1/2} \Big) \\
   & + \frac{1}{n^2} \sum_{g=1}^n \sum_{h \neq g} \sum_{\ell \neq g,h} \sum_{k \neq g,h,\ell}  a_g a_h b_\ell b_k \Big] 
  \end{align*}
  We apply the Cauchy-Schwarz inequality to $\sum_{\ell=1}^n a_\ell$ (and similarly for the $b_k$'s)
  which is bounded above by $( n \sum_{\ell=1}^n a_\ell^2 )^{1/2}$.
  Therefore, the terms in the previous display are bounded above by
  \begin{align}
    \frac{K}{n^2} \Big[ 
    \Big( \sum_{\ell=1}^n a_\ell^2 \Big) \Big( \sum_{k=1}^n b_k^2 \Big)
	+ \Big( \Big(\sum_{\ell=1}^n a_\ell^2\Big)^{3/2} \Big(\sum_{k=1}^n b_k^2\Big)^{1/2}
	+ \Big(\sum_{\ell=1}^n a_\ell^2\Big)^{1/2} \Big(\sum_{k=1}^n b_k^2\Big)^{3/2}\Big) \Big]
    \label{eq:tightness_bound} 
  \end{align}
  which we abbreviate by $K / n^2 ( AB + A^{3/2} B^{1/2} + A^{1/2} B^{3/2} )$.
  Now, verifying the unconditional Billingsley criterion implies conditional tightness in probability:
  Suppose $\vartheta_n^*$ is a random element of $D[s,\tau]$ for which tightness in probability given some random element $\vartheta_n$ shall be shown
  and for which it hold that:
  For all $\eta, \varepsilon > 0$ we have for the modulus of continuity $w''$, using the notation of  \cite{billingsley99}, that 
  $$ \pr (w''(\vartheta_n^*, \delta) \geq \varepsilon)  \leq \eta. $$
  Then, by the Markov inequality,
  \begin{align*}
   \pr ( \pr(w''(\vartheta_n^*, \delta) \geq \varepsilon \mid \vartheta_n ) \geq \varepsilon ) 
      & \leq E( \pr( w''(\vartheta_n^*, \delta) \geq \varepsilon \mid \vartheta_n )) /  \varepsilon \\
      & \leq \pr( w''(\vartheta_n^*, \delta) \geq \varepsilon ) /  \varepsilon
	\leq \eta / \varepsilon.
  \end{align*}

  Therefore, we need to find a suitable upper bound for the expectation of~\eqref{eq:tightness_bound}.
  By the Cauchy-Schwarz inequality we have the upper bound
  \begin{align*}
    \frac{K}{n^2} \Big[ E(AB) + \sqrt{E(A^2 + B^2) E(AB)} \Big].
  \end{align*}
  Here, we begin by looking for a bound of
  \begin{align*}
   E(AB) = \sum_{\ell \neq k} E(a_\ell^2 b_k^2) + \sum_{\ell=1}^n E(a_\ell^2 b_\ell^2) 
    = n(n-1) E(a_1^2) E(b_1^2) + n E(a_1^2 b_1^2) 
  \end{align*}
  We note that each $a_\ell$ and $b_\ell$ are differences of indicator functions.
  In particular, we have that
  \begin{align*}
   a_\ell^2 = & (1\{ X_\ell(s) \in \I, X_\ell(t) \in \J, C_\ell > t\} - 1\{ X_\ell(s) \in \I, X_\ell(t_1) \in \J, C_\ell > t_1\})^2 \\
     = & 1\{ X_\ell(s) \in \I, X_\ell(t) \in \J, t \geq C_\ell > t_1\} \\
     & + 1\{ X_\ell(s) \in \I, X_\ell(t_1) \notin \J, X_\ell(t) \in \J, C_\ell > t\} \\
     & + 1\{ X_\ell(s) \in \I, X_\ell(t_1) \J, X_\ell(t) \notin \J, C_\ell > t\}.
  \end{align*}
  Therefore, we can bound
  \begin{align*}
   a_\ell^2 \leq  1\{ t \geq C_\ell > t_1 \} + 1\{ X_\ell(t_1) \notin \J, X_\ell(t) \in \J \} + 1\{ X_\ell(t_1) \in \J, X_\ell(t) \notin \J \}.
  \end{align*}
  A similar bound holds for $b_\ell^2$. Thus, we have for the product
  \begin{align*}
   a_\ell^2 b_\ell^2 = & 1\{X_\ell(t) \notin \J, X_\ell(t_2) \in \J, t \geq  C_\ell > t_1\} \\
    + & 1\{X_\ell(t) \in \J, X_\ell(t_2) \notin \J, t \geq C_\ell > t_1\} \\
    + & 1\{X_\ell(t_1) \notin \J, X_\ell(t) \in \J, t_2 \geq C_\ell > t\} \\
    + & 1\{X_\ell(t_1) \in \J, X_\ell(t) \notin \J, t_2 \geq C_\ell > t\} \\
    + & 1\{X_\ell(t_1) \notin \J, X_\ell(t) \in \J, X_\ell(t_2) \notin \J \} \\
    + & 1\{X_\ell(t_1) \in \J, X_\ell(t) \notin \J, X_\ell(t_2) \in \J \}.
  \end{align*}
  Similarly to the proof in Appendix~\ref{sec:proofL1}, we obtain the upper bound
  \begin{align*}
   E(AB) \leq Kn^2 (H(t_2) - H(t_1))^2
  \end{align*}
  due to conditions (1) to (4) in the main article, where the function $H$ was defined in Appendix~\ref{sec:proofL1}.

  Next, we find an upper bound for 
  \begin{align*}
   E(A^2) & = E \Big[ \Big( \sum_{\ell=1}^n a_\ell^2 \Big)^2 \Big] = n(n-1) E(a_1^2)^2 + n E(a_1^2) \leq n^2 E(a_1^2) \\ 
   & \leq Kn^2 (H(t_2) - H(t_1)),
  \end{align*}
  which holds for the same reasons as before. 
  Similarly, we obtain the same upper bound for $E(B^2)$.
  We finally conclude that
  \begin{align*}
   E & \Big( |W^*_{\JI}(t) - W^*_{\JI}(t_1)|^2 |W^*_{\JI}(t_2) - W^*_{\JI}(t)|^2 \Big) \leq K (H(t_2) - H(t_1))^{3/2}.
  \end{align*}
  As argued above, Theorem~13.5 in \cite{billingsley99} now implies conditional tightness in probability given $\chi$.

  The same holds true for $W_\I^*$ which, as was the case for the original estimator, is more easily shown  than for $W_{\JI}^*$ and it is thus left to the reader.
  All in all, conditional convergence in distribution of $\sqrt{n}(P_{\I\J}^*(s,\cdot) - \wh P_{\I\J}(s,\cdot))$ on $D[s,\tau]$ given $\chi$ in probability follows by the functional delta-method.
  The Gaussian limit process is the same as for the original estimator.

\section{Proof of Theorem 3}  

  Consider subsequences of increasing sample sizes for which $n_1 / (n_1 + n_2) \rightarrow \lambda \in (0,1)$. 
  Due to the independence of the eventual limit distribution on the specific value of $\lambda$, 
  the convergence in distribution still holds under the bounded $\liminf$-$\limsup$-condition.
  The proof follows from a two-fold application of the continuous mapping theorem to the integration and the subtraction functional as well as from the independence of both samples.
  If $U$ is a zero-mean Gaussian process with covariance function $\Gamma_{\IJ}$,
  then $\int_s^\tau U(v) \d v$ has a zero-mean normal distribution with variance
  $ \int_s^\tau \int_s^\tau \Gamma_\IJ(u,v) \d u \d v. $

\section{Detailed Derivation of the Covariance Function $\Gamma_{\I\J}$}
\label{sec:cov}

As stated in Theorem~1, the covariance function $\Gamma_{\I\J}$ consists of several simpler covariance functions of which $cov(L_0(r), G_\ji(t))$ and $cov(L_1(r), G_\ji(t))$ still need to be determined.
Note that, if $L_2$ denotes the limit Gaussian process of the normalized Aalen--Johansen estimator for the cumulative incidence function $F_2$ of the second risk, then 
$$cov(L_0(r), G_\ji(t)) =  - cov(L_1(r), G_\ji(t)) - cov(L_2(r), G_\ji(t))$$
due to $F_0 + F_1 + F_2 = 1$ on $[s,\tau]$.
It is therefore enough to focus on the covariance function which involves the Aalen--Johansen limit.
This may be derived based on the asymptotic linear representation of all three separate estimators
  as Hadamard-derivatives $\phi'$, $\varphi'$, and $\psi'$, which are continuous and linear functionals, applied to an empirical process $\mb P_n$ and to $(\wh p_{\JI}, \wh p_\I)$, respectively.
  More precisely, we have
  \begin{align*}
   \sqrt{n}( \wh F_0 - F_0) & = \sqrt{n} \phi' ( \mb P_n - \mb P ) + o_p(1), \\
   \sqrt{n}( \wh F_1 - F_1) & = \sqrt{n} \varphi' ( \mb P_n - \mb P ) + o_p(1), \\
   \sqrt{n}( \wh p_{\ji} - p_{\ji}) & = \sqrt{n} \psi' ( \wh p_{\JI} - p_{\JI}, \wh p_\I - p_\I) + o_p(1).
  \end{align*}
  Here $\mb P_n$ is the empirical process based on the censoring or event times 
  $T_\ell = \inf\{t \geq s: Y_\ell(t) = 0\}$,
  the censoring or competing risk indicator $\delta_\ell = 1\{ T_{\ell1} < T_{\ell2} \} + 2 \cdot 1\{ T_{\ell1} > T_{\ell2} \}$,
  and $\varepsilon_\ell = 1\{ X_\ell(s) \in \I \}$ as well as $\eta_\ell = 1\{ Y_\ell(s) = 1 \}$,
  where
   $T_{\ell1} = \inf\{t \geq s: Z_\ell(t) = 1, Y_\ell(t) = 1\}$
   and $T_{\ell2} = \inf\{t \geq s: Z_\ell(t) = 2, Y_\ell(t) = 1 \}$.
   Finally, $\mb P = P^{(T_\ell, \delta_\ell, \varepsilon_\ell, \eta_\ell)}$ is the distribution of each single observed quadruple.
   The empirical process is indexed by 
   $f_{1z}(t, d, x, y) = 1\{ t \leq z, d = 1, x = 1, y = 1 \}$, $f_{2z}(t, d, x, y) = 1\{ t \leq z, d = 2, x = 1, y = 1 \}$, and $g_{z}(t, d, x, y) = 1\{ t \geq z, x = 1, y = 1 \}$,
   $z \geq s$.
   This certainly is a Donsker class; see Example~3.9.19 in \cite{vaart96} for a similar indexing leading to the (all-cause) Nelson--Aalen estimator.
   In the following part, we calculate the exact structure of the above linear expansions for the case $n=1$.
   This suffices to determine the asymptotic covariance structure of the estimators
   as the covariances of the above linear functionals are the same for all $n$.
   
   As argued in the main article, it does not matter whether the left- or right-continuous versions of the at-risk indicators $Y_\ell$ are considered.
   This is also true for the at risk functions involved in the Kaplan--Meier and Aalen--Johansen estimators.
   For technical convenience, we choose the right-continuous at-risk functions.
   Let $p_{\mc A_\J, \mc R_\J}(t) = P(T_1 \leq t, \delta_1 \in \{1,2\}, X_1(s) \in \I, Y_1(s) = 1) $,
   $p_{\mc A_\J}(t) = P(T_1 \leq t, \delta_1 = 1, X_1(s) \in \I, Y_1(s) = 1) $,
   and $p_{\mc R_\J}(t) = P(T_1 \leq t, \delta_1 = 2, X_1(s) \in \I, Y_1(s) = 1) $.
   The above Hadamard-derivatives are given by
   \begin{align*}
    (\phi'(\mb P))(t) = ( \phi_{1, - A_0}' \circ \phi'_{2, p_{\mc A_\J, \mc R_\J}, 1/p_\I} \circ \phi'_{3, p_{\mc A_\J, \mc R_\J}, p_\I} )(\mb P)(t).
   \end{align*}
   Thereof, $\phi_{1, - A_0}'$ is the Hadamard-derivative of the product-integral leading to the Kaplan--Meier estimator, i.e.
   \begin{align*}
    (\phi_{1, - A_0}' (h))(t) = - \prodi_{s < u \leq t} (1 - A_0(\d u)) \int_s^t \frac{\d h(u)}{1 - \Delta A_0(u)}
    = - F_0(t) (h(t) - h(s));
   \end{align*}
   cf. e.g. Lemma~3.9.30 in \cite{vaart96}, but take notice of the additional minus sign.
   Here, $A_0$ is the continuous cumulative (all-cause) hazard function for all individuals for which $X(s) \in \I$.
   The functional $\phi'_{2, p_{\mc A_\J, \mc R_\J}, 1/p_\I}$ is the Hadamard-derivative of the Wilcoxon functional, i.e.
   \begin{align*}
    (\phi'_{2, p_{\mc A_\J, \mc R_\J}, 1/p_\I} (h_1, h_2))(t) = \int_s^t h_2(u) \d p_{\mc A_\J, \mc R_\J}(u)
     + \int_s^t \frac{1}{p_\I(u)} \d h_1(u);
   \end{align*}
   cf. e.g. Example~3.9.19 in \cite{vaart96}. Finally, $\phi'_{3, p_{\mc A_\J, \mc R_\J}, p_\I}$ is the Hadamard-derivative of $(a,b) \mapsto (a, 1/b)$, after indexing the distribution with $(f_{1 \cdot} + f_{2 \cdot}, g_\cdot)$, i.e.
   $$ (\phi'_{3, p_{\mc A_\J, \mc R_\J}, p_\I}(\mb P))(t) = \Big( \mb P (f_{1 t} + f_{2 t}), -\frac{\mb P g_t}{p_\I^2(t)} \Big).  $$
   From the above derivations, we have that
   \begin{align*}
    (\phi'(\mb P_1 & - \mb P))(t_1) =  F_0(t_1) \Big[ \int_s^{t_1} \frac{1\{ Y(u) = 1, X(s) \in \I \} - p_\I(u) }{p^2_\I(u)} \d p_{\mc A_\J, \mc R_\J}(u) \\
     & - \int_s^{t_1} \frac{1}{p_\I(u)} ( 1\{ X(\d u) \in \mc A_\J \cup \mc R_\J, \delta \neq 0, X(s) \in \I \} - p_{\mc A_\J, \mc R_\J}(\d u)) \Big] \\
     & = F_0(t_1) \Big[ \int_s^{t_1} \frac{1\{ Y(u) = 1, X(s) \in \I \}}{p^2_\I(u)} \d p_{\mc A_\J, \mc R_\J}(u) \\
     & \quad  - \int_s^{t_1} \frac{1\{ X(\d u) \in \mc A_\J \cup \mc R_\J, \delta \neq 0, X(s) \in \I \}}{p_\I(u)} \Big] \\
     & = \varepsilon \eta F_0(t_1) \Big[ \int_s^{T \wedge t_1} \frac{\d p_{\mc A_\J, \mc R_\J}(u)}{p^2_\I(u)}
      - \frac{1\{ T \leq t_1, \delta \neq 0 \}}{p_\I(T)} \Big]
   \end{align*}

   We proceed similarly for the Aalen--Johansen estimator 
   $$\wh F_1(t) = \int_s^t \wh F_0(u-) \d \wh A_1(u)$$
   of $F_1(t) = \int_s^t F_0(u-) \d A_1(u)$, where now $A_1$ denotes the first cause-specific cumulative hazard of those individuals for which $X(s) \in \I$.
   Define $p_{\mc A_\J}(t) = P(T \leq t, \delta = 1, X(s) \in \I, Y(s) = 1) $.
   Thus, as the estimator $\wh F_1$ consists of the Wilcoxon functional applied to the left-continuous version of a product integral and another Wilcoxon functional (the Nelson--Aalen estimator), we have
   \begin{align*}
    (\varphi'(\mb P_1 - \mb P))(t_2) = (\tilde \phi'_{2, A_1, F_0} \circ (\phi'_{2, p_{\mc A_\J}, 1/p_\I} \circ \tilde \phi'_{3, p_{\mc A_\J}, p_\I}, \phi'))(\mb P_1 - \mb P)(t_2).
   \end{align*}
    The tilde on top of the first Hadamard-derivative shall express the slight variation of the original functional, as we now use the left-continuous version of the integrand.
    Furthermore,
    $\tilde \phi'_{3, p_{\mc A_\J}, p_\I}(\mb P_1 - \mb P)(t)
    = (\mb P f_{1t} - p_{\mc A_\J}(t), - \frac{\mb P g_t - p_\I(t)}{p_\I^2(t)})$.
    Hence, the previous display can be expanded to
   \begin{align*}
    & \int_s^{t_2} \phi'(\mb P_1 - \mb P) (u-) \d A_1(u) + \int_s^{t_2} F_0(u-) \d (\phi'_{2, p_{\mc A_\J}, 1/p_\I} \circ \tilde \phi'_{3, p_{\mc A_\J}, p_\I})(\mb P_1 - \mb P)(u) \\
    & = \int_s^{t_2} F_0(u-) \Big[ \int_s^{u-} \frac{1\{ Y(v) = 1, X(s) \in \I \}}{p^2_\I(v)} \d p_{\mc A_\J, \mc R_\J}(v) \\
    & \quad - \int_s^{u-} \frac{1\{ X(\d v) \in \mc A_\J \cup \mc R_\J, \delta \neq 0, X(s) \in \I \}}{p_\I(v)} \Big] \d A_1(u) \\
    & \quad + \int_s^{t_2} F_0(u-) \d \Big[ - \int_s^{\cdot} \frac{1\{ Y(v) = 1, X(s) \in \I \}}{p^2_\I(v)} \d p_{\mc A_\J}(v) \\
    & \quad + \int_s^{\cdot} \frac{1\{ X(\d v) \in \mc A_\J, \delta \neq 0, X(s) \in \I \}}{p_\I(v)} \Big] \\
    & = \varepsilon \eta \Big( \int_s^{t_2} \Big[ \int_s^{(T \wedge u)-} \frac{\d p_{\mc A_\J, \mc R_\J}(v)}{p^2_\I(v)} - \frac{1\{T < u, \delta \neq 0\}}{p_\I(T)} \Big] \d F_1(u) \\
    & \quad - \int_s^{T \wedge t_2} \frac{F_0(u-)}{p^2_\I(u)} \d p_{\mc A_\J}(u) 
    + 1\{T \leq t_2, \delta = 1 \} \frac{ F_0(T-) }{p_\I(T)} \Big) \\
    & = \varepsilon \eta \Big( F_1(t_2) \int_s^{T \wedge t_2} \frac{\d p_{\mc A_\J, \mc R_\J}(u)}{p^2_\I(u)} 
    - \int_s^{T \wedge t_2} \frac{ F_1(u) }{p_\I^2(u)} \d p_{\mc A_\J, \mc R_\J}(u) \\ 
    & \quad - \frac{1\{T \leq t_2, \delta \neq 0\}}{p_\I(T)} ( F_1(t_2) - F_1(T) ) \\
    & \quad - \int_s^{T \wedge t_2} \frac{F_0(u)}{p^2_\I(u)} \d p_{\mc A_\J}(u) 
    + 1\{ T \leq t_2, \delta = 1 \} \frac{ F_0(T) }{p_\I(T)} \Big) \\
    & = \varepsilon \eta \Big( \int_s^{T \wedge t_2} \frac{F_1(t_2) - 1 + F_2(u)}{p^2_\I(u)} \d p_{\mc A_\J, \mc R_\J}(u)
    + \int_s^{T \wedge t_2} \frac{F_0(u)}{p^2_\I(u)} \d p_{\mc R_\J}(u) \\
    & \quad - \frac{1\{T \leq t_2, \delta \neq 0\}}{p_\I(T)} ( F_1(t_2) - F_1(T) ) 
    + 1\{ T \leq t_2, \delta = 1 \} \frac{ F_0(T) }{p_\I(T)} \Big)
   \end{align*}
   
   Finally, we calculate the remaining linear functional applied to a pair of centered indicator functions that appear in the estimator $\wh p_\ji$:
   \begin{align*}
    & (\psi'(1\{ X(s) \in \I, Y(s) = 1, X(\cdot) \in \J, Y(\cdot)  = 1 \}, \\
    & \qquad 1\{ X(s) \in \I, Y(s) = 1, Y(\cdot)  = 1 \}))(t_3) \\
    & = \frac{1\{ X(s) \in \I, Y(s) = 1, X(t_3) \in \J, Y(t_3)  = 1 \} - p_\JI(t_3)}{p_\I(t_3)} \\
    & \quad  - \frac{p_\JI(t_3)}{p_\I^2(t_3)} (1\{ X(s) \in \I, Y(s) = 1, Y(t_3)  = 1 \} - p_\I(t_3)) \\
    & = p_\I^{-1}(t_3) ( 1\{ X(s) \in \I, Y(s) = 1, X(t_3) \in \J, Y(t_3)  = 1 \} \\
    & \quad - p_\ji(t_3) 1\{ X(s) \in \I, Y(s) = 1, Y(t_3)  = 1 \} ) \\
    & = \frac{\varepsilon \eta \ 1\{ T > t_3 \} }{p_\I(t_3)} ( 1\{ X(t_3) \in \J \} - p_\ji(t_3))   
   \end{align*}

   For later calculations we note that
   $$ P(T \in \d t , \delta = 1, X(s) \in \I, Y(s) = 1) = S_C(t) F_1(\d t) P(X(s)\in \I), $$
   that $P^{(T,\varepsilon, \eta)}(\d t, 1, 1) = - p_\I(\d t) $, and that
   \begin{align*}
    &P(T \in \d t , \delta = 1, X(t_3) \in \J, X(s) \in \I, Y(s) = 1) \\
   &= S_C(t) P_{\I\J1}(s,t_3,\d t) P(X(s)\in \I),
   \end{align*}
   where $P_{\I\J1}(s,t_3, t) = P(X(t) = 1, X(t_3) \in \J \mid X(s) \in \I)$.
   
   After this preparation, we are now able to determine the asymptotic covariance.
   Because of the independencies and the identical distribution structure in the asymptotic, linear representation,
   \begin{align*}
   & cov (L_1(t_2),  G_{\ji}(t_3)) \\
     & = n \cdot cov((\varphi'(\mb P_n - \mb P))(t_2),  (\psi'(\wh p_\JI - p_\JI, \wh p_\I - p_\I))(t_3)) \\
     & =  cov((\varphi'(\mb P_1 - \mb P))(t_2), (\psi'(1\{ X(s) \in \I, Y(s) = 1, X(\cdot) \in \J, Y(\cdot)  = 1 \} - p_{\JI}, \\
     & \qquad 1\{ X(s) \in \I, Y(s) = 1, Y(\cdot)  = 1 \} - p_\I))(t_3)) \\
    & = E \Big(  \varepsilon \eta \frac{1\{ T > t_3 \} }{p_\I(t_3)}
    \Big[ 
      \int_s^{T \wedge t_2} \frac{F_1(t_2) - 1 + F_2(u)}{p^2_\I(u)} \d p_{\mc A_\J, \mc R_\J}(u) \\
    & \quad + \int_s^{T \wedge t_2} \frac{F_0(u)}{p^2_\I(u)} \d p_{\mc R_\J}(u)
    - \frac{1\{T \leq t_2, \delta \neq 0\}}{p_\I(T)} ( F_1(t_2) - F_1(T) )  \\
    & \quad + 1\{ T \leq t_2, \delta = 1 \} \frac{ F_0(T) }{p_\I(T)} 
    \Big]
    \Big[ 
       1\{ X(t_3) \in \J \} - p_\ji(t_3)
    \Big]
    \Big)
   \end{align*}
    We begin by considering the case $t_2 \leq t_3$ for which the covariance reduces to
     \begin{align*}
     & E \Big(  \varepsilon \eta \frac{1\{ T > t_3 \} }{p_\I(t_3)}
    \Big[ 
      \int_s^{t_2} \frac{F_1(t_2) - 1 + F_2(u)}{p^2_\I(u)} \d p_{\mc A_\J, \mc R_\J}(u)
    + \int_s^{t_2} \frac{F_0(u)}{p^2_\I(u)} \d p_{\mc R_\J}(u) 
    \Big] \\
    & \quad \times \Big[ 
       1\{ X(t_3) \in \J \} - p_\ji(t_3)
    \Big]
    \Big) \\
    & = \frac{1}{p_\I(t_3)} \Big[ 
      \int_s^{t_2} \frac{F_1(t_2) - 1 + F_2(u)}{p^2_\I(u)} \d p_{\mc A_\J, \mc R_\J}(u)
    + \int_s^{t_2} \frac{F_0(u)}{p^2_\I(u)} \d p_{\mc R_\J}(u) 
    \Big] \\
    & \quad \times \Big( 
      \int_{t_3}^\infty P^{(T,1\{X(t_3) \in \J\}, \varepsilon, \eta )}(\d v, 1, 1, 1)
       - p_\ji(t_3) \int_{t_3}^\infty  P^{(T, \varepsilon, \eta )}(\d v, 1, 1) 
    \Big).
    \end{align*}
    The round bracket in the last line equals $ p_\JI(t_3) - p_\ji(t_3)  p_\I(t_3) = 0$.
    Hence, asymptotic independence follows in case of $t_2 \leq t_3$.
    
    We hence continue with the case $t_2 > t_3$.
    Let $f$ and $g$ be any suitable functions.
    The following side calculation due to a decomposition will prove beneficial:
    \begin{align*}
     & E \Big(1\{ T > t_3 \} \int_s^{T \wedge t_2} f(u) \d g(u) \Big)  \\
     & \quad = \int_{t_3}^{t_2} \int_s^v f(u) \d g(u) \d P^T(v) 
      + \int_{t_2}^\infty \int_s^{t_2} f(u) \d g(u) \d P^T(v) 
    \end{align*}
    For a better presentation of the calculations, we split the covariance terms into two parts.
    First,
   \begin{align*}
     & E \Big(  \varepsilon \eta \frac{1\{ T > t_3 \} }{p_\I(t_3)}
    \Big[ 
      \int_s^{T \wedge t_2} \frac{F_1(t_2) - 1 + F_2(u)}{p^2_\I(u)} \d p_{\mc A_\J, \mc R_\J}(u)
    + \int_s^{T \wedge t_2} \frac{F_0(u)}{p^2_\I(u)} \d p_{\mc R_\J}(u) 
    \Big] \\
    & \quad \times \Big[ 
       1\{ X(t_3) \in \J \} - p_\ji(t_3)
    \Big]
    \Big) \\
    & = \frac{1}{p_\I(t_3)} \Big( 
      \int_{t_3}^{t_2} \int_s^{v} \frac{F_1(t_2) - 1 + F_2(u)}{p^2_\I(u)} \d p_{\mc A_\J, \mc R_\J}(u)  P^{(T,1\{X(t_3) \in \J\}, \varepsilon, \eta )}(\d v, 1, 1, 1) \\
      & \quad + \int_{t_3}^{t_2} \int_s^{v} \frac{F_0(u)}{p^2_\I(u)} \d p_{\mc R_\J}(u)  P^{(T,1\{X(t_3) \in \J\}, \varepsilon, \eta )}(\d v, 1, 1, 1) \\
      & \quad - p_\ji(t_3) \int_{t_3}^{t_2} \int_s^{v} \frac{F_1(t_2) - 1 + F_2(u)}{p^2_\I(u)} \d p_{\mc A_\J, \mc R_\J}(u)  P^{(T, \varepsilon, \eta )}(\d v, 1, 1) \\
      & \quad - p_\ji(t_3) \int_{t_3}^{t_2} \int_s^{v} \frac{F_0(u)}{p^2_\I(u)} \d p_{\mc R_\J}(u)  P^{(T,\varepsilon, \eta )}(\d v, 1, 1) \\
      & \quad + \int_{t_2}^\infty \int_s^{t_2} \frac{F_1(t_2) - 1 + F_2(u)}{p^2_\I(u)} \d p_{\mc A_\J, \mc R_\J}(u)  P^{(T,1\{X(t_3) \in \J\}, \varepsilon, \eta )}(\d v, 1, 1, 1) \\
      & \quad + \int_{t_2}^\infty \int_s^{t_2} \frac{F_0(u)}{p^2_\I(u)} \d p_{\mc R_\J}(u)  P^{(T,1\{X(t_3) \in \J\}, \varepsilon, \eta )}(\d v, 1, 1, 1) \\
      & \quad - p_\ji(t_3) \int_{t_2}^\infty \int_s^{t_2} \frac{F_1(t_2) - 1 + F_2(u)}{p^2_\I(u)} \d p_{\mc A_\J, \mc R_\J}(u)  P^{(T, \varepsilon, \eta )}(\d v, 1, 1) \\
      & \quad - p_\ji(t_3) \int_{t_2}^\infty \int_s^{t_2} \frac{F_0(u)}{p^2_\I(u)} \d p_{\mc R_\J}(u)  P^{(T,\varepsilon, \eta )}(\d v, 1, 1)
    \Big).
   \end{align*}
   Thereof, the sum of the last four terms can be summarized similarly as in the case $t_2 \leq t_3$, i.e. it equals
   \begin{align*}
    & \frac{1}{p_\I(t_3)} \Big[ 
      \int_s^{t_2} \frac{F_1(t_2) - 1 + F_2(u)}{p^2_\I(u)} \d p_{\mc A_\J, \mc R_\J}(u)
    + \int_s^{t_2} \frac{F_0(u)}{p^2_\I(u)} \d p_{\mc R_\J}(u) 
    \Big] \\
    & \quad \times \Big( S_C(t_2) P(X(s) \in \I) P_{\I\J0}(s,t_3,t_2) - p_\ji(t_3) F_0(t_2) S_C(t_2) P(X(s) \in \I) \Big) \\
    & = \frac{S_C(t_2)}{p_\I(t_3)} \Big[ 
     - \int_s^{t_2} \frac{1}{S_C(u)} \frac{F_1(t_2) - 1 + F_2(u)}{F_0^2(u)} \d F_0(u)
    + \int_s^{t_2} \frac{1}{S_C(u)} \frac{\d F_2(u)}{F_0(u)} 
    \Big] \\
    & \quad \times \Big( P_{\I\J0}(s,t_3,t_2) - p_\ji(t_3) F_0(t_2) \Big).
   \end{align*}
    The sum of the remaining four terms equals, after integration by parts,
  \begin{align*}
    & \frac{1}{p_\I(t_3)} \Big( 
      - \int_s^{t_2} \frac{F_1(t_2) - 1 + F_2(u)}{p^2_\I(u)} \d p_{\mc A_\J, \mc R_\J}(u) P_{\I\J0}(s,t_3,t_2) S_C(t_2) P(X(s)\in \I) \\
      & \quad + \int_s^{t_3} \frac{F_1(t_2) - 1 + F_2(u)}{p^2_\I(u)} \d p_{\mc A_\J, \mc R_\J}(u) P_{\I\J0}(s,t_3,t_3) S_C(t_3) P(X(s)\in \I) \\
      & \quad + \int_{t_3}^{t_2} \frac{F_1(t_2) - 1 + F_2(u)}{p^2_\I(u)} P_{\I\J0}(s,t_3,u) S_C(u) P(X(s)\in \I) \d p_{\mc A_\J, \mc R_\J}(u) \\
      & - \int_s^{t_2} \frac{F_0(u)}{p^2_\I(u)} \d p_{\mc R_\J}(u) P_{\I\J0}(s,t_3,t_2) S_C(t_2) P(X(s)\in \I) \\
      & \quad + \int_s^{t_3} \frac{F_0(u)}{p^2_\I(u)} \d p_{\mc R_\J}(u)  P_{\I\J0}(s,t_3,t_3) S_C(t_3) P(X(s)\in \I)  \\
      & \quad + \int_{t_3}^{t_2} \frac{F_0(u)}{p^2_\I(u)}  P_{\I\J0}(s,t_3,u) S_C(u) P(X(s)\in \I) \d p_{\mc R_\J}(u) \\
      & \quad + p_\ji(t_3) \int_s^{t_2} \frac{F_1(t_2) - 1 + F_2(u)}{p^2_\I(u)} \d p_{\mc A_\J, \mc R_\J}(u) p_\I(t_2) \\
      & \quad - p_\ji(t_3) \int_s^{t_3} \frac{F_1(t_2) - 1 + F_2(u)}{p^2_\I(u)} \d p_{\mc A_\J, \mc R_\J}(u) p_\I(t_3) \\
      & \quad - p_\ji(t_3) \int_{t_3}^{t_2} \frac{F_1(t_2) - 1 + F_2(u)}{p_\I(u)} \d p_{\mc A_\J, \mc R_\J}(u) \\
      & \quad + p_\ji(t_3) \int_s^{t_2} \frac{F_0(u)}{p^2_\I(u)} \d p_{\mc R_\J}(u) p_\I(t_2)
	- p_\ji(t_3) \int_s^{t_3} \frac{F_0(u)}{p^2_\I(u)} \d p_{\mc R_\J}(u) p_\I(t_3) \\ 
	& \quad - p_\ji(t_3) \int_{t_3}^{t_2} \frac{F_0(u)}{p_\I(u)} \d p_{\mc R_\J}(u)
    \Big) \\
    & = \frac{1}{p_\I(t_3)} \Big( 
       \int_s^{t_2} \frac{S_C(t_2)}{S_C(u)} \frac{F_1(t_2) - 1 + F_2(u)}{F_0^2(u)} \d F_0(u) P_{\I\J0}(s,t_3,t_2) \\
      & \quad - \int_s^{t_3} \frac{S_C(t_3)}{S_C(u)} \frac{F_1(t_2) - 1 + F_2(u)}{F_0^2(u)} \d F_0(u) P_{\I\J0}(s,t_3,t_3) \\
     & \quad  - \int_{t_3}^{t_2} \frac{F_1(t_2) - 1 + F_2(u)}{F_0^2(u)} P_{\I\J0}(s,t_3,u) \d F_0(u) \\
      & \quad - \int_s^{t_2} \frac{S_C(t_2)}{S_C(u)} \frac{1}{F_0(u)} \d F_2(u) P_{\I\J0}(s,t_3,t_2) 
       + \int_s^{t_3} \frac{S_C(t_3)}{S_C(u)} \frac{\d F_2(u)}{F_0(u)}  P_{\I\J0}(s,t_3,t_3) \\  
      & \quad  + \int_{t_3}^{t_2} \frac{P_{\I\J0}(s,t_3,u)}{F_0(u)}  \d F_2(u) \\
      & \quad 
	- p_\ji(t_3)  \int_s^{t_2}  \frac{S_C(t_2)}{S_C(u)} \frac{F_1(t_2) - 1 + F_2(u)}{F_0^2(u)}  F_0(\d u) F_0(t_2) \\
      & \quad + p_\ji(t_3) \int_s^{t_3} \frac{S_C(t_3)}{S_C(u)} \frac{F_1(t_2) - 1 + F_2(u)}{F_0^2(u)}  F_0(\d u) F_0(t_3) \\
      & \quad  + p_\ji(t_3) \int_{t_3}^{t_2} \frac{F_1(t_2) - 1 + F_2(u)}{F_0(u)} \d F_0(u) \\
      & \quad + p_\ji(t_3) \int_s^{t_2} \frac{S_C(t_2)}{S_C(u)} \frac{\d F_2(u)}{F_0(u)}  F_0(t_2)
	- p_\ji(t_3) \int_s^{t_3} \frac{S_C(t_3)}{S_C(u)} \frac{\d F_2(u)}{F_0(u)} F_0(t_3) \\
      & \quad - p_\ji(t_3) (F_2(t_2) - F_2(t_3))
    \Big) 
   \end{align*}
   Comparing these terms with those four which we derived first, we see that all integrals from $s$ to $t_2$ are cancelled out.
   Therefore, the following terms remain:
   \begin{align*}
    & \frac{1}{p_\I(t_3)} \Big( 
      - \int_s^{t_3} \frac{S_C(t_3)}{S_C(u)} \frac{F_1(t_2) - 1 + F_2(u)}{F_0^2(u)} \d F_0(u) P_{\I\J0}(s,t_3,t_3) \\
      & \quad - \int_{t_3}^{t_2} \frac{F_1(t_2) - 1 + F_2(u)}{F_0^2(u)} P_{\I\J0}(s,t_3,u) \d F_0(u) \\
      & \quad + \int_s^{t_3} \frac{S_C(t_3)}{S_C(u)} \frac{\d F_2(u)}{F_0(u)}  P_{\I\J0}(s,t_3,t_3)  
       + \int_{t_3}^{t_2} \frac{P_{\I\J0}(s,t_3,u)}{F_0(u)}  \d F_2(u) \\
      & \quad + p_\ji(t_3) \int_s^{t_3} \frac{S_C(t_3)}{S_C(u)} \frac{F_1(t_2) - 1 + F_2(u)}{F_0^2(u)}  F_0(\d u) F_0(t_3) \\
      & \quad + p_\ji(t_3) \int_{t_3}^{t_2} \frac{F_1(t_2) - 1 + F_2(u)}{F_0(u)} \d F_0(u) \\
      & \quad - p_\ji(t_3) \int_s^{t_3} \frac{S_C(t_3)}{S_C(u)} \frac{\d F_2(u)}{F_0(u)} F_0(t_3) 
	- p_\ji(t_3) (F_2(t_2) - F_2(t_3))
    \Big) 
   \end{align*}
   Due to $P_{\I\J0}(s,t_3,t_3) = p_{\ji}(t_3) F_0(t_3)$, this equals
   \begin{align*}
        & \frac{1}{p_\I(t_3)} \Big( 
      - \int_{t_3}^{t_2} \frac{F_1(t_2) - 1 + F_2(u)}{F_0^2(u)} P_{\I\J0}(s,t_3,u) \d F_0(u) 
       + \int_{t_3}^{t_2} \frac{P_{\I\J0}(s,t_3,u)}{F_0(u)}  \d F_2(u) \\
      & \quad  + p_\ji(t_3) \int_{t_3}^{t_2} \frac{F_1(t_2) - 1 + F_2(u)}{F_0(u)} \d F_0(u) 
	- p_\ji(t_3) (F_2(t_2) - F_2(t_3))
	    \Big) = A
   \end{align*}
  We continue by deriving the second half of the asymptotic covariance function also in case of $t_3 < t_2$, i.e.
   \begin{align*}
    &  E \Big(  \varepsilon \eta \frac{1\{ T > t_3 \} }{p_\I(t_3)}
    \Big[  - \frac{1\{T \leq t_2, \delta \neq 0\}}{p_\I(T)} ( F_1(t_2) - F_1(T) ) 
    + 1\{ T \leq t_2, \delta = 1 \} \frac{ F_0(T) }{p_\I(T)} 
    \Big] \\
    & \quad \times \Big[ 
       1\{ X(t_3) \in \J \} - p_\ji(t_3)
    \Big]
    \Big) \\
    & = \frac{1}{p_\I(t_3)} \Big( 
    - \int_{t_3}^{t_2} \frac{F_1(t_2) - F_1(u)}{p_\I(u)} P^{(T, \delta, 1\{X(t_3) \in \J\}, \varepsilon, \eta)}(\d u , \{1,2\}, 1, 1, 1) \\
    & \quad + p_\ji(t_3) \int_{t_3}^{t_2} \frac{F_1(t_2) - F_1(u)}{p_\I(u)} P^{(T, \delta, \varepsilon, \eta)}(\d u , \{1,2\}, 1, 1) \\
    & \quad + \int_{t_3}^{t_2} \frac{F_0(u)}{p_\I(u)} P^{(T, \delta, 1\{X(t_3) \in \J\}, \varepsilon, \eta)}(\d u ,1,1, 1, 1) \\
    & \quad - p_\ji(t_3) \int_{t_3}^{t_2} \frac{F_0(u)}{p_\I(u)} P^{(T, \delta, \varepsilon, \eta)}(\d u , 1, 1, 1) \\
    & = \frac{1}{p_\I(t_3)} \Big( 
    \int_{t_3}^{t_2} \frac{F_1(t_2) - F_1(u)}{F_0(u)} P_{\I\J0}(s,t_3,\d u) 
     - p_\ji(t_3) \int_{t_3}^{t_2} \frac{F_1(t_2) - F_1(u)}{F_0(u)} \d F_0(u) \\
    & \quad + P_{\I\J1}(s,t_3,t_2) - P_{\I\J1}(s,t_3,t_3)
     - p_\ji(t_3) (F_1(t_2) - F_1(t_3))  \Big) = B
   \end{align*}
   We apply integration by parts to the first integral. Note that, if $u$ is the variable of integration,
   \begin{align*}
    \d \frac{F_1(t_2) - F_1(u)}{F_0(u)}
    & = - \frac{F_1(t_2) \d F_0(u)}{F_0^2(u)}
      - \frac{\d F_1(u)}{F_0(u)} 
      + \frac{F_1(u) \d F_0(u)}{F_0^2(u)} \\
     &  = \frac{(- F_1(t_2) + F_1(u)) \d F_0(u)}{F_0^2(u)}
      + \frac{\d F_0(u)}{F_0(u)} 
      + \frac{\d F_2(u)}{F_0(u)} 
   \end{align*}
    Therefore, the first integral of $B$ equals
    \begin{align*}
    &  - (F_1(t_2) - F_1(t_3)) p_\ji(t_3)
     + \int_{t_3}^{t_2} P_{\I\J0}(s,t_3,u) \frac{F_1(t_2) - F_1(u)}{F_0^2(u)} \d F_0(u) \\
    & - \int_{t_3}^{t_2} P_{\I\J0}(s,t_3,u) \frac{\d F_0(u)}{F_0(u)} 
     - \int_{t_3}^{t_2} P_{\I\J0}(s,t_3,u) \frac{\d F_2(u)}{F_0(u)}
    \end{align*}
    All in all, the sum of both parts $A$ and $B$ equals
    \begin{align*}
     & \frac{1}{p_\I(t_3)} \Big( 
      \int_{t_3}^{t_2} \frac{F_0(u)}{F_0^2(u)} P_{\I\J0}(s,t_3,u) \d F_0(u) 
       - p_\ji(t_3) \int_{t_3}^{t_2} \frac{F_0(u)}{F_0(u)} \d F_0(u) \\
    & \quad + p_\ji(t_3) (F_0(t_2) - F_0(t_3)) 
      - \int_{t_3}^{t_2} P_{\I\J0}(s,t_3,u) \frac{\d F_0(u)}{F_0(u)} 
      + P_{\I\J1}(s,t_3,t_2) \\
     & \quad - P_{\I\J1}(s,t_3,t_3)
      - p_\ji(t_3) (F_1(t_2) - F_1(t_3))\Big) \\
     & =  \frac{1}{p_\I(t_3)} \Big(
       P_{\I\J1}(s,t_3,t_2) - P_{\I\J1}(s,t_3,t_3) 
      - p_\ji(t_3) (F_1(t_2) - F_1(t_3))\Big),
    \end{align*}
    which, due to $ P_{\I\J 1}(s,t_3,t_3) = F_1(t_3) $, further simplifies to
   \begin{align*}
      \frac{1\{t_3 < t_2 \}}{p_\I(t_3)} \Big(
       P_{\I\J1}(s,t_3,t_2) - F_1(t_3) - p_\ji(t_3) (F_1(t_2) - F_1(t_3))\Big).
   \end{align*}
{\blue
\section{Conditions (1)--(4) in the Simulation Study}
In this section we prove that the conditions (1)--(4) of the main manuscript are satisfied in the situation of Section~4 therein.
Recall that the stochastic development of the process $X$ depends on the state occupied at time $s=4$. Hence, conditional on $X(4)=1$ and after time $s=5$, $X$ is a Markov process with cumulative hazard matrix
$$ A(t-5) = \begin{pmatrix}
             -0.32 & 0.3 & 0.02 \\
             0.3 & -0.4 & 0.1 \\
             0 & 0 & 0 
            \end{pmatrix} \cdot (t-5), \quad t \geq s.
 $$
 Otherwise, if $X(4)=0$, the process $X$ subsequently is a Markov process with cumulative hazard matrix
$$ B(t-5) = \begin{pmatrix}
             -0.62 & 0.6 & 0.02 \\
             0.3 & -0.4 & 0.1 \\
             0 & 0 & 0 
            \end{pmatrix} \cdot (t-5), \quad t \geq s.
 $$
 Write $p= P(X(4) = 1)$ and $q = 1-p$.
 In this conditionally homogeneous Markov set-up it follows that the matrix $P(5,t)$ of transition probabilities at time $t\geq 5$  is given by
 $$ P(5,t) = p \cdot \exp( A(t-5) ) + q \cdot \exp(B(t-5)), $$
 where $\exp(D) = \sum_{k=0}^\infty \frac1{k!}D^k$ is the matrix exponential of a square matrix $D$.
 Now, it is obvious that $t \mapsto P(5,t)$ is infinitely often continuously differentiable, hence (1) and (2) are satisfied.
 In the same way, (3) and (4) are satisfied as well because $X$ is a Markov process conditionally on $X(4)$, so its state at an earlier intermediate time $t_1 > s = 5 > 4$ is irrelevant.}

\end{document}